\documentclass[12pt]{article}
\usepackage[margin=1in]{geometry} 

\usepackage[utf8]{inputenc}
\usepackage{lipsum}

\usepackage{amssymb,amsmath,bm}
\usepackage[colorlinks=true, urlcolor=blue,linkcolor=blue, citecolor=blue]{hyperref}

\usepackage{amsthm}
\usepackage{enumerate}
\usepackage{enumitem}   
\usepackage{appendix}
\usepackage{romannum}

\author{Mine Çağlar\\
Department of Mathematics\\ 
Ko\c{c} University\\
Istanbul, Turkey \\
Email: mcaglar@ku.edu.tr
\and 
İhsan Demirel\\
Department of Mathematics\\ 
Ko\c{c} University\\
Istanbul, Turkey \\
E-mail: idemirel17@ku.edu.tr}

\title{Backward Monge Potential and   Monge-Amp\'ere Equation}

\begin{document}
	\maketitle

\begin{abstract}
In this paper, Monge-Kantorovich problem is considered  in  the  infinite  dimension on an abstract Wiener space  $(W,H,\mu)$, where  $H$ is  Cameron-Martin space and $\mu$ is the Gaussian measure.  We  study the regularity of optimal transport maps  with a quadratic cost function assuming that both initial and target measures have a strictly positive Radon-Nikodym density  with  respect  to  $\mu$. Under conditions on the density functions, the forward and backward transport maps can be written in terms of  Sobolev derivative of so-called Monge-Brenier maps, or Monge potentials.   We show  Sobolev regularity of the backward potential  under the assumption that  the density of the initial measure is log-concave and prove   that it  solves Monge-Ampére equation.
\end{abstract}

\renewcommand{\thefootnote}{}



\renewcommand{\thefootnote}{\arabic{footnote}}
\setcounter{footnote}{0}

	\pagenumbering{arabic}

	\newcommand{\nn}{\nonumber}
	\newcommand{\To}{\rightarrow}
	\newcommand{\supp}{\mbox{supp}}

	\newcommand{\myred}[1]{\textcolor{red}{#1}}

	\newcommand{\nnumber}{\stepcounter{equation}\tag{\theequation}\label}
	\newcommand{\vertiii}[1]{{\left\vert\kern-0.25ex\left\vert\kern-0.25ex\left\vert #1 
			\right\vert\kern-0.25ex\right\vert\kern-0.25ex\right\vert}}

	\newcommand{\R}{\mathbb{R}}
	\newcommand{\N}{\mathbb{N}}
	\newcommand{\Z}{\mathbb{Z}}
	\newcommand{\Q}{\mathbb{Q}}
	\newcommand{\D}{\mathbb{D}}
	\newcommand{\E}{\mathbb{E}}
	
	\newcommand{\FF}{\mathcal{F}}
	\newcommand{\LL}{\mathcal{L}}
	\newcommand{\XX}{\mathcal{X}}		
	\newcommand{\oname}[1]{\operatorname{#1}}

	\newtheorem{theorem}{Theorem}[section]
	\newtheorem{corollary}{Corollary}[section]
	\newtheorem{conjecture}{Conjecture}[section]
	\newtheorem{proposition}{Proposition}[section]
	\newtheorem{lemma}{Lemma}[section]
	\newtheorem{definition}{Definition}[section]
	\newtheorem{example}{Example}[section]
	\newtheorem{remark}{Remark}[section]
	\newtheorem{solution}{Solution}[section]
	\newtheorem{case}{Case}[section]
	\newtheorem{problem}{Problem}[section]
	\newtheorem{condition}{Condition}[section]

	\section{Introduction}
	
  Monge problem is motivated by the application of moving given piles of sand to fill up  holes of the same  volume  with a minimum total cost of transportation \cite{ambrosio_2003,villani2003topics}.  Piles and holes are modeled by  probability measures $\rho$ and $\nu$ defined on some measurable sets $X$ and $Y$, respectively, with a measurable and non-negative cost function   $c:X\times Y\rightarrow \R_+ \cup\{+ \infty \}.$ The aim is to find a transport map $T:X\rightarrow Y$, with  $T\rho=\nu$, that minimizes the expected cost of moving the sand. A solution to Monge problem may not always exist, as it does not allow split of mass.  Monge-Kantorovich problem overcomes this difficulty by the relaxation that mass from location $x\in X$ can be moved possibly to several locations in $Y$ with the objective of finding a transference plan $\gamma^*\in\Gamma(\rho,\nu) $ that minimizes the functional
		$$J[\gamma]=\int_{X\times Y} C(x,y) d\gamma(x,y) $$
where $\Gamma(\mu,\nu)$ is the set of joint probability measures on $X\times Y$ with first and second marginals   $\rho$ and $\nu$, respectively.	If $X$ and $Y$ are Polish spaces and   $c$ is lower semi-continuous, then	Monge-Kantorovich problem admits a minimizer \cite[Thm.1.5]{ambrosio91}. In the infinite dimension,   Monge-Kantorovich problem  has been studied on Wiener space first in  \cite{Feyel_Ustunel_2004} with a quadratic cost function. Later in \cite{Cavalletti2012,fang,nolot2013}, different cost functions are considered  with  other initial and target measures. 

In this paper, we consider an abstract Wiener space $(W,H,\mu)$, where $H$ is  Cameron-Martin space and $\mu$ is the Gaussian measure, and  define	$\text{d}\rho=e^{-f}\text{d}\mu$ and $\text{d}\nu=e^{-g}\text{d}\mu$ for measurable $f,g:W\to\mathbb{R}$. Under conditions on $f$ and $g$,  there exists  so-called forward Monge potential or Monge-Brenier map $\varphi $ such that	 $T = I_W+\nabla\varphi$ is the optimal transport map, where $\nabla  $ denotes the Gross-Sobolev derivative operator. Moreover, the inverse $S$ of $T$ is given	by $S=I_W+\nabla\psi$ for some $\psi: W\to\mathbb{R}$, called backward Monge potential, satisfying $(S\times I_W)\rho = (I_W\times T)\nu = \gamma$. Sobolev regularity of the Monge-Brenier maps is an important issue in order to write  Jacobian functions associated with transformations $T$ and $S$.    In Theorem \ref{4.1}, assuming $f$ is $(1-c)$-convex, equivalently $e^{-f}$ is log-concave, we show that  $\nabla^2\psi$ is in $ L^2(\nu,H\otimes H)$ and provide  an upper estimate for its norm as 
	\begin{align*}
		\E_{\nu}\left[|\nabla^2\psi|_2^2 \right]
		&\leq \frac{3}{c}\left( \E_{\rho}\left[|\nabla\varphi|_H^2 \right] +\E_{\nu}\left[|\nabla g|_H^2\right]+\E_{\rho}\left[|\nabla f|_H^2\right]  \right)\, ,
		\end{align*}
where $H\otimes H$ denotes the space of Hilbert-Schmidt operators on $H$ and $|\cdot|_2$ denotes the Hilbert-Schmidt norm. We use the approximation approach of \cite{ustunel2017variational}, where the initial measure $\rho$ is taken as $\mu$ and the target measure $\nu$ is assumed to be absolutely continuous with respect to $\mu$.   

  Monge potentials are closely related to Monge-Amp\'ere equation defined in Wiener space. On $\R^n$, with  $F,G:\R^n\rightarrow \R_+$, Monge-Amp\`ere equation in $T$ can  be written  as 
	\begin{align*}
	F=G \circ T \det J_T \nnumber{fine_monge_ampere}
	\end{align*}
where $J_T$ is the Jacobian of $T$. When  $F$ and $G$ are the densities of $\rho$ and $\nu$ with respect to Lebesgue measure, respectively, the corresponding transport map $T$ is a solution. In finite dimension, regularity of solution of Monge-Kantorovich problem and Monge-Amp\'ere equation have been studied in \cite{Caffarelli1992,Figalli2013,Philippis2013,MCCANN1997153}. For the aim of defining an analogous equation in Wiener space, suppose  $F=e^{-f}$ and $G=e^{-g}$. In other words, Radon-Nikodym derivatives of $\rho$ ve $\nu$ with respect to Lebesgue measure are proportional to $e^{-f-\frac{|x|^2}{2}}$ and $e^{-g-\frac{|x|^2}{2}}$, respectively,  for $x\in\R^n$ and \eqref{fine_monge_ampere} becomes
	\begin{align*}
	e^{-f(x)-\frac{|x|^2}{2}}=e^{-g\circ T(x)-\frac{|T(x)|^2}{2}}{\det} J_T(x) .
	\end{align*}
On $\R^n$, we have  $T = I_W+J_\varphi$, where $J_\varphi$ denotes the Jacobian of  forward Monge potential $\varphi$.
Using the identities $J_ T= I_{\R^n}+\bigtriangleup \varphi$ and
	$\det(I_{\R^n}+\cdot)={\det}_2(I_{\R^n}+\cdot) e^{\text{tr}(\cdot)}$, 
	where ${\det}_2$ denotes Carleman-Fredholm determinant, tr denotes the trace of a matrix and $\bigtriangleup$ denotes the Laplacian, we get
	\begin{align*}
	e^{-f}=e^{-g\circ T}{\det}_2 (I_{\R^n}+\nabla^2 \varphi) e^{\text{tr}(\nabla^2 \varphi)}e^{-x J_\varphi} e^{-\frac{|J_\varphi|^2}{2}} .
	\end{align*}
In view of the identity 
$-\mathcal{L}\varphi=\Delta \varphi-x J_\varphi$ $=\text{tr}(\nabla^2 \varphi)-x J_\varphi$,
where $\mathcal{L}$ is the generator of the Ornstein–Uhlenbeck process on $\R^n$,  Monge-Amp\'ere equation can now be written as
\begin{align*}
	e^{-f}=e^{-g\circ T}{\det}_2\left( I_{\R^n} +\bigtriangleup \varphi\right) \exp \left[ -\mathcal{L}\varphi -\frac{1}{2}|J_\varphi|^2 \right] .
	\end{align*}
Since only ${\det}_2$ is well-defined in the infinite dimension, this form can be used to define Monge-Amp\'ere equation in Wiener space by replacement of $|\cdot|$ with the norm in $H$, $J_\varphi$ with Gross-Sobolev derivative $\nabla \varphi$, and $\bigtriangleup \varphi$ with $\nabla^2\varphi $. In \cite{Feyel_Ustunel_2004}, it is proved that  Monge potential $\varphi $ solves Monge-Amp\'ere equation under the condition $f$ and $g$ are two positive random variables with values
in a bounded interval $[a, b]$.  Monge-Amp\'ere equation for $\varphi$ has been obtained in \cite{FEYEL2004,FEYEL2006}  with $f=0$ and $g$  an $H$-convex Wiener
function. In \cite{KOLESNIKO2005,bogachev2013}, the authors have shown that $\varphi$ satisfies  
Monge-Amp\'ere equation for  $g\in L^1(
\mu)$ and $  ge^{-g} \in L^1(\mu
)$ as well as the case  when $e^{-f/2}\in \D_{2,1}$ and $g=0$. The case when both initial and target measures  are absolutely continuous with respect to $\mu$ is studied in \cite{fang} for bounded $f$ and lower bounded $g$ with the additional assumption that  second Sobolev derivative of $g$ exists. 

In the settings of the present paper, the backward potential is shown to be regular in Theorem \ref{4.1}, that is, $\nabla^2 \psi \in L^2(\nu,H\otimes H)$ under the assumption that the initial measure $\rho$ has a log-concave density with respect to Gaussian measure $\mu$. Therefore, we prove in Theorem \ref{4.2}   that backward potential $\psi$  solves Monge-Amp\'ere equation given by
	\begin{align*}
		e^{-g}=e^{-f\circ S}{\det}_2\left( I_{H} +\nabla^2\psi\right) \exp \left[ -\mathcal{L}\psi -\frac{1}{2}|\nabla\psi|^2_H \right] 
		\end{align*}
$\nu$-almost surely.

The organization of the paper is as follows. In Section 2, we review the preliminary definitions and introduce the notation used in the paper. Section 3 focuses on the finite dimensional and smooth case as a basis for the infinite dimension. In Section 4, the regularity of backward potential is shown and Monge-Amp\'ere equation is considered. 
	
	\section{Preliminaries}
	
	Let $(W,H,\mu)$  be an abstract Wiener  space. The corresponding Cameron-Martin space is denoted by $H$. The norm in $H$ will be denoted by $|\cdot|_H$. 
	If $h \in H$,
	then there exists $(l_n) \subset W^*$ such that image of this sequence under injection $W^*\hookrightarrow H$, say $(\tilde{l}_n )$, converges to $h$
	in $H$. Therefore, the sequence of random variables  $(\langle l_n,\cdot \rangle_H) $ is Cauchy in $L^p(\mu)$ for any $p\geq0$. We denote its limit by $\delta h$, which is $N(0,|h|^2_H)$ random variable.

	The function $F:W\to \R$ is called a cylindrical
	Wiener functional if it is of the form
	\begin{align*}
	F(\omega)=f\left(\delta{h_1}(\omega), \ldots, \delta{h_n}(\omega)\right), \quad h_{1} \cdots h_{n} \in H, \quad f \in \mathcal{S}\left(\mathbb{R}^{n}\right)
	\end{align*}
	for some $n \in \mathbb{N}$ and $\mathcal{S}\left(\mathbb{R}^{n}\right)$ denotes the Schwartz space of rapidly decreasing functions on $\mathbb{R}^{n}$  \cite{Feyel_Ustunel_2004,ustunel1995}.  We denote the collection of Wiener functionals by $\mathcal{S}(W )$. For such $F\in \mathcal{S}(W )$ and $h \in H,$ we define
	$$
	\nabla_{h} F(\omega)=\left.\frac{d}{d \epsilon} F(\omega+\epsilon h)\right|_{\epsilon=0}
	$$
	For fixed $\omega\in W$, $h \mapsto \nabla_{h}F(\omega)$ is continuous and linear on $H$. Therefore, there exists an element in $H$, say $\nabla F$ such that $\nabla_hF=\langle\nabla F, h\rangle_H.$ The operator $\nabla:F\mapsto \nabla F$
	is linear  from $\mathcal{S}(W)$ into the space of
	$H$-valued Wiener functionals $L^p(\mu; H)$ for any $p > 1.$
	The operator $\nabla$ is closable from $L^p(\mu)$ into $L^p(\mu;H)$ for any $p>1.$ The completion of $S(W)$ under the norm
	\begin{align*}
	\| \cdot\|_{p,1}= \|\cdot\|_{L^p(\mu)}+ \|\cdot\|_{L^p(\mu;H)} 
	\end{align*}  
	is denoted by $\D_{p,1}$, which is a Banach space with norm $\| \cdot\|_{p,1}$ \cite{ustunel1995}.
	The definition can be also extended to the collection of Wiener functionals which take values in a separable Hilbert space $\mathcal{X}.$ The completion  of the collection of $\mathcal{X}$-valued Wiener functionals $\mathcal{S}(W; \mathcal{X})$ under norm 
	\begin{align*}
	\|\cdot\|_{L^p(\mu;\mathcal{X})}+ \|\cdot\|_{L^p(\mu;\mathcal{X}\otimes H)} 
	\end{align*}  
	is denoted by $\D_{p,1}(\mathcal{X}).$ Higher order derivatives can also be defined, e.g. we say $F\in \D_{p,2}$ if $\nabla F\in \D_{p,1}(H)$ and write $\nabla^2F=\nabla(\nabla F).$
	
	Let $\nu$ be measure on $W$ absolutely continuous with respect to $\mu$ with Radon-Nykodm  derivative $L.$ If $\int_W |\nabla L|^2_He^{-L}\ d\mu<\infty,$ then operator $\nabla$ is closable over $\mathcal{S}(W)$ under the norm 
	\begin{align*}
	\|\cdot\|_{L^p(\nu)}+ \|\cdot\|_{L^p(\nu; H)}\ .
	\end{align*}
	We will denote the completion by $\D_{p,1}(\nu).$ 
	The adjoint of continuous and
	linear operator $\nabla$ will be denoted by $\delta$. That is, for suitable $\xi:W\to H$ and any $F\in\D_{p,1}$, we have
	$$ \E\left[\langle\nabla F , \xi\rangle _{H}\right]=\E[\varphi \cdot \delta \xi].$$	
	The operator $\delta: \xi\mapsto \delta\xi$ is called the divergence operator. We can also define $\delta_\nu$ by this procedure, i.e, $\delta_\nu$ is the adjoint of
	the Sobolev derivative $\nabla$ under the  measure $\nu$.
	
	For measurable $f :W\to\R$ and $t\geq 0$ the Ornstein-Uhlenbeck semi-group $\left(P_{t}\right)_{t \geq 0} $   is given   by 
	$$
	\left(P_{t} f\right)(x)=\int_{W} f\left(e^{-t} x+\sqrt{1-e^{-2 t}} y\right) d \mu(y).
	$$
	Its infinitesimal generator
	is denoted by $-\LL$ and we call $\LL$ the Ornstein–Uhlenbeck operator.  The norm given by $\| (I+\LL)^{r/2}(\cdot)\|_{L^p(\mu)}
	$
	is equivalent to the norm $\|\cdot\|_{p,r}$ for any $p>1$ and $r\in \N.$ Completion of $\mathcal{S} ((W)$ with respect to the norm $\| (I+\LL)^{r/2}(\cdot)\|_{L^p(\mu)}
	$ again denoted by $\D_{p,r}$ which has $\D_{q,-r},$ $q^{-1}=1-p^{-1}$, as its continuous dual for any $p>1,$ $r\in\R.$ Similarly, we can define  $\D_{p,r}(\XX)$ as a completion of $\mathcal{S} (W; \XX)$ with respect to the norm
	$\| (I+\LL)^{k/2}F\|_{L^p(\mu;\XX)}$ and its continuous dual is $\D_{q,-r}(\XX^*),$ where $\XX^*$ is the dual of $\XX.$ The space 
	\begin{align*}
	\D(\XX)=\bigcap_{p,r} \D_{p,r}(\XX)
	\end{align*}
	is dense in $\D_{p,r}(\XX),$ for any $p>1$ and $r\in\R.$ If the sequence $(F_n)$ converges to zero in $\D_{p,r}(\XX),$ for any $p>1$ and $r\in\R,$ we say that it converges to zero in $\D(\XX).$ Under this topology $\D(\XX)$ is a complete, locally convex topological vector
	space. The continuous dual $\D^*(\XX^*)$ of $\D(\XX)$ is called Meyer-Watanabe distributions on $W$ with values in $\XX^*.$
	
	A measurable function $f: W \rightarrow \mathbb{R} \cup\{\infty\}$ is called $\alpha$-convex, $\alpha\in \R$,  if the map
	$$
	h \mapsto  f(\cdot+h)+\frac{\alpha}{2}|h|_{H}^{2}
	$$
	is convex on the Cameron-Martin space $H$ with values in $L^{0}(\mu) $ a.s.

 For Monge-Kantorovich problem in the Wiener space, the cost function $C(x,y):W\times W\to \R^+\cup\{+\infty\}$ is taken as 
	$$
	C(x,y)=\begin{cases} 
	|x-y|_H^2 & \text{if }x-y\in H \\
	+\infty& \text{otherwise}.
	\end{cases}
	$$ 
 in \cite{Feyel_Ustunel_2004}.
The following result is the starting point of the present work \cite[Thm.1.1]{fang}.
	\begin{theorem}\label{existence_MKP}
		Let $\rho$, $\nu$ be the probability measures on $W$ such that
$		\text{d}\rho=F \text{d}\mu $  and $ \text{d}\nu=G  \text{d}\mu$,
		where $\mu$ is the Wiener measure, $F:W\to\mathbb{R}$ and $G:W\to\mathbb{R}$ are measurable functions such that 
		\begin{align*}
		\int_{W}\frac{|  \nabla F|_H^2}{F}\;{d}\mu<\infty \text{ and } \int_{W}\frac{|  \nabla G|_H^2}{G}\;{d}\mu<\infty 
		\end{align*}
		and the function $F $ satisfies Poincar{\'e} inequality, that is, for every cylindrical functional $\xi :W\to \mathbb{R}$
		\begin{align*}
		(1-c)\int_{W}( \xi - E_\rho [\xi])^2F\;{d}\mu\leq \int_{W}|\nabla \xi|_H^2F\;{d}\mu \nnumber{poincare}
		\end{align*}
		for some $c\in [0,1).$
Then, there exists a $\varphi \in \mathbb{D}_{2,1}(\rho)$ such that $ T=I_W+\nabla \varphi$ is the unique solution of Monge problem for $(\rho, \nu)$ and the probability measure $\gamma$ given by
		$
		\gamma=(I_W \times T)\rho 
		$
		is the unique solution of Monge-Kantorovich problem for $(\rho, \nu).$
		Moreover, $T$ is $\rho$-a.s invertible and the inverse map $S=T^{-1}$ has the form  $S=I_W+\nabla\psi$, 
		where $\psi \in \mathbb{D}_{2,1}(\nu)$.
	\end{theorem}
	Finally, let $\mathcal{P}(X )$ denote  the probability measures on a measurable space $X$.  For  $m_1, m_2 \in \mathcal{P}(X )$ on a Polish  space $({X},d)$, the Wasserstein distance of order $p \in[1, \infty) $ between $m_1$ and $m_2$ is defined by 
	\begin{align*}
	d_{p}(m_1, m_2)=\left(\inf _{\gamma \in \Gamma(m_1, m_2)} \int_{{X}} d(x, y)^{p}\ d \gamma(x, y)\right)^{1 / p}. \nnumber{Wasserstein distance}
	\end{align*}
	Note that $d_p$ is not a metric, since it can take the value of $+\infty $. However, its restriction on a subset
	of $\mathcal{P}(X ) \times \mathcal{P}(X )$ where it takes finite values leads to a a complete metric space. 
	The relative entropy of $m_1$ with respect to $m_2$ given by
	\begin{align*}
	H(m_1\mid m_2)=\begin{cases} 
	\int_X \frac{dm_1}{dm_2}\log(\frac{dm_1}{dm_2}) \ dm_2 & \text{if $m_1\ll m_2$  } \\
	+\infty& \text{otherwise}.
	\end{cases}\nnumber{relative entropy}
	\end{align*}
	We take $X=W$ and $d(x,y)=|x-y|_H$ and $p=2.$  The following form of Talagrand's inequality, derived in \cite[Thm.3.1]{Feyel_Ustunel_2004},  provides a nice relation between  Wasserstein distance and the relative entropy. For any probability measure $\mu_0$ on $W$, it holds that
	\begin{align*}
	d_2^2(\mu_0,\mu)\leq 2H(\mu_0 \mid \mu) \nnumber{Talagrand}.
	\end{align*}
	
	\section{An Estimate for Smooth Backward Monge Potential}
	
	In this section,  we assume that  forward Monge potential $\varphi$ is smooth, that is,  
	$\varphi \in \D_{2,k}$ for all $k\in\N$. We first prove an identity  for forward Monge potential $\varphi$ and then find an upper bound for $\D_{2,2}(\nu)$-norm of $\psi$, which will be auxiliary for regularity in the infinite dimension as shown in the next section. Similar results  have been obtained in \cite{ustunel2017variational} when the  initial measure is Gaussian and the target measure is absolutely continuous with respect to Gaussian measure.  
	
	\begin{remark}\label{smooth_condition}  When  $f$ and $g$ are smooth and bounded from below, and $\mu$ is a standard Gaussian measure on $\R^n$, then  $\varphi$ is smooth. This would follow from \cite{alex2010sobolev} and \cite[Thm.4.14 ]{villani2003topics}.
	\end{remark}
	\begin{proposition}\label{smooth_relation}
		When $W$ is finite dimensional and $\varphi$ is smooth,  Monge potential $\varphi$ satisfies the  relation
		\begin{align*}
		\nabla \varphi+\nabla f\circ T-\nabla g=\delta_\rho\left[\left(I_{H}+\nabla^2\varphi \right)^{-1} -I_{H}  \right] \; .
		\end{align*}
	\end{proposition}
	\begin{proof}
		From \cite[Prop.2.5 ]{var_rep}, it is known that for bounded $ f,g$
		\begin{align*}
		-\log\int_{W}e^{-g+f}\;{d}\rho=\inf\left( \int_{W}(g-f)\ {d}m + H(\gamma|\rho):\; m\in \mathcal{P}(W)\right) 
		\end{align*}
		The infimum is attained at $\nu$ given that $H\left( \nu|\rho\right) <\infty.$ Moreover, since $T\rho=\nu$, we have
		\begin{align*}
		-\log\int_{W}e^{-g+f}\;{d}\rho
		&=\inf\left( \int_{W}(g-f)\circ (I_W+\nabla a)\ {d}\rho + H\left( (I_W+\nabla a) \rho\mid \rho\right):\ a\in\mathbb{D}_{2,1}(\rho)\right) 
		\\
		&\geq \inf\left( \int_{W}(g-f)\circ (I_W+ \xi)\ {d}\rho + H\left( (I_W+ \xi)\rho\mid \rho\right) :\ \xi\in\mathbb{D}_{2,0}(\rho)\right) 
		\\
		&\geq\inf\left( \int_{W}(g-f)\ {d}m + H(m\mid\rho):\ m\in \mathcal{P}(W)\right) .
		\end{align*}
Therefore,  $\inf\left\{ K_{g-f}(\xi):\ \xi\in\mathbb{D}_{2,0}(\rho) \right\} $ is attained at $\xi=\nabla\varphi,$ where 
		\begin{align*}
		K_{g-f}(\xi)= \int_{W}(g-f)\circ (I_W+\xi)\;{d}\rho + H\left( (I_W+\xi)\rho\mid \rho\right) . \nnumber{K_maps}
		\end{align*}
We have $\nu\ll\rho$ with $\ d\nu/d\rho=e^{-g}/e^{-f}.$ 
It follows that
		\begin{align*}
		H\left( T \rho\mid\rho\right)
		&=H(\nu|\rho)\\
		&=\int_{W}e^{-g}(-g+f))\;d\mu\\
		&=\int_{W}e^{-g\circ T}(-g+f)\circ T \ \Lambda_{\varphi} \;d\mu \nnumber{relative_ent}
		\end{align*}	where $\Lambda_{\varphi}$ is the Gaussian Jacobian of $T=I_W+{\nabla\varphi}$ given by
		\begin{align*}
		\Lambda_{T}=\operatorname{det}_2(I_{H}+{\nabla^2\varphi})\exp\left( -\mathcal{L}\varphi-\frac{1}{2}|{\nabla\varphi}|_H^2\right)  .
		\end{align*}
	Since we have smooth $\varphi\in\mathbb{D}_{2,1}(\rho)$, it solves Monge-Amp{\'e}re equation, that is,
		\begin{align*}
		e^{-f}=e^{-g\circ T }\operatorname{det}_2(I_{H}+{\nabla^2\varphi})\exp\left( -\mathcal{L}\varphi-\frac{1}{2}|{\nabla\varphi}|_H^2\right) 
		\end{align*}
		If we substitute this to \eqref{relative_ent}, we get 
		\begin{align*}
		H\left( T \rho\mid\rho\right)
		&=\int_{W}e^{-f}(-f-\log \Lambda_{\varphi} +f\circ T )\ d\mu
		\end{align*}
		Consider the map $T_t=I_W+t\nabla\varphi$ for $t\in[0,1).$ Here, $T_t$ is strongly $(1-t)-$monotone shift \cite[Lem.6.2.1]{Ustunel_Zakai_2000}. Moreover, for $\xi\in \mathbb{D}_{2,1}(H)$ with $||\nabla\xi||_2$ is in the space of bounded random variable $L^\infty(\rho)$  and sufficiently small $\epsilon>0$, the shift $T_{t,\epsilon}=I_W+t\nabla\varphi+\epsilon\xi$ is still strongly monotone. Therefore, $\Lambda_{t\nabla\varphi+\epsilon\xi}>0$ holds a.s. \cite{Ustunel_Zakai_2000}. Moreover, $T_{t,\epsilon}\rho$ is absolutely continuous with respect to $\rho$ \cite[Thm.7.3]{Feyel_Ustunel_2004} and we have 
		\begin{align*}
		e^{-f}=e^{-g_{t,\epsilon}\circ T_{t,\epsilon} }\Lambda_{T_{t,\epsilon}} 
		\end{align*}
		where $e^{g_{t,\epsilon}}={dT_{t,\epsilon}\rho}{/}{d\rho}$ $\cite{MCCANN1997153}$.	
		A similar calculation shows that 
		\begin{align*}
		H\left( T_{t,\epsilon} \rho|\rho\right)
		&=\int_{W}e^{-f}(-f+f\circ T_{t,\epsilon}-\log \Lambda_{T_{t,\epsilon}})\;d\mu\nnumber{entropy_t_epsilon}
		\end{align*}
		and
		\begin{align*}
		K_{g_t-f}(t\nabla\varphi+\epsilon\xi)= \int_{W}g_t\circ T_{t,\epsilon}   -f-\log \Lambda_{T_{t,\epsilon}} \;{d}\rho, \nnumber{J_maps1}
		\end{align*} where $e^{g_{t}}={dT_{t}\rho}/{d\rho}.$
		Note that $t\varphi$ is unique (up to a constant) Monge potential of   MKP for $\Sigma(\rho,T_{t}\rho)$ and it minimizes  $K_{f_t-g}$  among all absolutely continuous shifts. Therefore, we have
		\begin{align*}
		\frac{d}{d\epsilon}K_{g_t-f}(t\nabla\varphi+\epsilon\xi)\Big|_{\epsilon=0}=0 .
		\end{align*}
		Observe that
		\begin{align*}
		\frac{d}{d\epsilon}\left( \int_W g \;{d}\rho\right)  \Big|_{\epsilon=0}
		&=0,		\\
		\frac{d}{d\epsilon}\left( \int_W g_{t}\circ T_{t,\epsilon} \;{d}\rho\right)  \Big|_{\epsilon=0}
		&=\int_{W}(\nabla g_t\circ T_t,\xi)\;{d}\rho\\
		\frac{d}{d\epsilon}\left( \int_W-\log \Lambda_{T_{t,\epsilon}}) \;{d}\rho\right)  \Big|_{\epsilon=0}
		&=\int_W-\operatorname{trace}\left(\left((I_{H}+t\nabla^2\varphi )^{-1} -I_{H} \right)\cdot \nabla\xi   \right) +\delta\xi+\left(\nabla\varphi,\xi \right)   \;{d}\rho\ ,
		\end{align*}
		for each $\xi\in\mathbb{D}_{2,1}(H)$ with $||\xi||_2\in L^\infty(\rho).$ Here, the first derivative is $0$ since it does not depend on $\epsilon,$ the second derivative result comes from definition and for the last derivative we have used \cite[Thm A.2.2]{Ustunel_Zakai_2000}.
		The set $$\Theta=\left\lbrace\xi\in\mathbb{D}_{2,1}(H): ||\xi||_2\in L^\infty(\rho) \right\rbrace $$ is dense in $L^p(\rho)$, we have 
		\begin{align*}
		t\nabla \varphi+\nabla g_t\circ T_t-\nabla f=\delta_\rho\left[\left(I_{H}+t\nabla^2\varphi \right)^{-1} -I_{ H}  \right] 
		\end{align*}
		Since $$e^{-f}=e^{-g_{t}\circ T_{t} }\Lambda_{t\nabla\varphi} $$
		$\ g_t\circ T_t$ converges in probability to $\nabla g\circ T$ as $t\to 1$, and we get  
		\begin{align*}
		\nabla \varphi+\nabla g\circ T-\nabla f=\delta_\rho\left[\left(I_{H}+\nabla^2\varphi \right)^{-1} -I_{H}  \right] .
		\end{align*}
	\end{proof}
	We need the following technical lemmas in order to prove the regularity of backward Monge potential. We refer to \cite{ustunel2017variational} for their proofs.
	\begin{lemma}{\label{positive_trace}}
	 Suppose $W$ is finite dimensional and $\varphi$ is smooth. Let $K=(I_H+\nabla^2\phi)^{-1}$ and $h\in H$, then we have 
		\begin{align*}
		\oname{trace}(K\nabla^3\varphi Kh\cdot K\nabla^3\varphi Kh)\geq0  \quad \quad \mu-a.s.
		\end{align*} 
	\end{lemma}
	\begin{lemma} If
		$\xi : W\to H$ is smooth, then
    $\delta_\rho\xi=\delta\xi+\langle\nabla f,\xi\rangle_H$.
	\end{lemma}
	\begin{lemma}If \label{smooth_exeptation}
		 $\xi : W\to H$ is smooth, then
		\begin{align*}
		\E_\rho[(\delta_\rho\xi)^2]=\E_\rho[(I_H+\nabla^2g,\xi\otimes\xi)_{H^{\otimes2}}+\oname{trace}(\nabla\xi\cdot\nabla\xi)].
		\end{align*}
	\end{lemma}
	\begin{proposition}
		\label{regulary_smooth_bakcward}
		Suppose $W$ is finite dimensional, $\varphi$ is smooth, and the function $f$ is $(1-c)$-convex  for some $c \in [0,1)$. Then, we have
		\begin{align*}
		c \E_{\nu}\left[|\nabla^2\psi|_2^2 \right]&\leq 3\left( \E_{\rho}\left[|\nabla\varphi|_H^2 \right] +\E_{\nu}\left[|\nabla g|_H^2\right]+\E_{\rho}\left[|\nabla f|_H^2\right]  \right)
		\end{align*}
	\end{proposition}
	\begin{proof}
		Taking expectation of the second moment of the norm of the expression in Proposition  \ref{smooth_relation} yields
		\begin{align*}
		\E_{\rho}\left[ \left| \delta_\rho\left(\left(I_{H}+\nabla^2\varphi \right)^{-1} -I_{H}  \right)\right|_H^2  \right]&\leq 3\left( \E_{\rho}\left[|\nabla\varphi|_H^2 \right] +\E_{\rho}\left[|\nabla g\circ T|_H^2\right]+\E_{\rho}\left[|\nabla f|_H^2\right]  \right)\\
		&=3\left( \E_{\rho}\left[|\nabla\varphi|_H^2 \right] +\E_{\nu}\left[|\nabla g|_H^2\right]+\E_{\rho}\left[|\nabla f|_H^2\right]  \right)\nnumber{bound1}
		\end{align*}
		Denote $\left(I_{H}+\nabla^2\varphi\right)^{-1}$  by $M.$ If we apply  Lemma \ref{smooth_exeptation}, then 
		\begin{align*}
		E_{\rho}\left[|\delta_{\rho}(M-I_{H})|_H^2 \right] &=\sum_{k=1}^{\infty}\E_{\rho}\left[\left( \delta_{\rho}(M-I_{H}) (e_k)\right)^2 \right]\\
		&=\sum_{k=1}^{\infty}\E_\rho\left[ \left( I_H+\nabla^2 f,(M-I_H) e_k\otimes (M-I_H) e_k \right)_{H{\otimes}H}  \right] \\ &{\,\,\,\,}+ \sum_{k=1}^{\infty}\E_{\rho}\left[ \text{trace}(\nabla(Me_k)\cdot\nabla(Me_k))\right]
		\end{align*}
		By Lemma \ref{positive_trace}, the second term in the last equality is positive. So, we have
		\begin{align*}
		\E_{\rho}\left[|\delta_{\rho}(M-I_{H})|_H^2 \right]
		&\geq\sum_{k=1}^{\infty}\E_{\rho}\left[ \left( I_H+\nabla^2f,(M-I_{H}) e_k\otimes (M-I_{H}) e_k \right)_{H{\otimes}H}  \right]\\
		&\geq c \sum_{k=1}^{\infty} \E_{\rho}\left[|(M-I_{H})e_k|_H^2 \right]\\
		&=c E_{\rho}\left[|(M-I_{H})|_2^2 \right]\ . \nnumber{bound2}
		\end{align*}  
		Since $T=I_W+\nabla\varphi$ and $S=I_W+\nabla\psi$ are inverses of each other, we have
		\begin{align*}
		(I_H+\nabla^2\varphi)^{-1}=(I_H+\nabla^2\psi)\circ T.\nnumber{inverse_relation1}
		\end{align*}
		Combining  \eqref{bound1},  \eqref{bound2} and \eqref{inverse_relation1}, we get 
		\begin{align*}
		c \E_{\nu}\left[|\nabla^2\psi|_2^2 \right]&\leq 3\left( \E_{\rho}\left[|\nabla\varphi|_H^2 \right] +\E_{\nu}\left[|\nabla g|_H^2\right]+\E_{\rho}\left[|\nabla f|_H^2\right]  \right).
		\end{align*}
	\end{proof}
	\section{Regularity  and Monge-Amp\`ere Equation}
	Recall that  $f:W\to\mathbb{R}$ and $g:W\to\mathbb{R}$ are measurable functions such that $f,g\in \D_{2,1}$ and 
	\begin{align}  \label{conds}
	\int_{W}|\nabla f|^2e^{-f}\;{d}\mu<\infty \quad \text{and} \quad  \int_{W}|\nabla g|^2e^{-g}\;{d}\mu<\infty
	\end{align}
	with $e^{-f} $ satisfying Poincar{\'e} inequality \eqref{poincare}, and the initial and target measures of Monge-Kantorovich problem are 	  $\text{d}\rho=e^{-f}\text{d}\mu$ and $\text{d}\nu=e^{-g}\text{d}\mu$.
	In this section, we will show that backward Monge potential $\psi$ solves Monge-Amp{\'e}re equation, that is, we have 
	\begin{align*}
	e^{-g}=e^{-f\circ S}{\det}_2\left( I_{H} +\nabla^2\psi\right) \exp \left[ -\mathcal{L}\psi -\frac{1}{2}|\nabla \psi|^2_H \right] \, 
	\end{align*}
	where  $S=T^{-1}=I_W+\nabla \psi$. 	However, we only know $\psi\in \D_{2,1}(\nu)$ so far. Therefore, we first establish an upper bound for $\D_{2,2}(\nu)$-norm of $\psi$, which implies that the second Sobolev derivative of $\psi$ is well-defined, in subsections \ref{s4.1} and \ref{s4.2}. Monge-Amp\'ere equation is considered in Subsection \ref{s4.3}.

	\subsection{Approximation Lemmas} \label{s4.1}
	
	In order to show the regularity of  backward Monge potential $\psi$ in the general setting of this section, we will approximate $f$ and $g$ of \eqref{conds}. This will be accomplished through several lemmas. 
	In Lemma \ref{lemma1}, we define a sequence of probability measures $\rho_n$ and $\nu_n$ which are absolutely continuous with respect to the standard Gaussian measure on $\mathbb{R}^n$, that is, in finite dimension. 
	
	\begin{lemma}\label{lemma1} Let
		$f,g\in\mathbb{D}_{2,1}$ satisfy \eqref{conds} with $e^{-f} $ satisfying Poincar{\'e} inequality \eqref{poincare}. Let $(\varphi,\psi)$ be the Monge potentials associated to the Monge-Kantorovitch problem $(\rho,\nu),$ where 
		$\text{d}\rho=e^{-f}\text{d}\mu \text{ and }  \text{d}\nu=e^{-g}\text{d}\mu.$
		Define $g_n$ and $f_n$ as  
		$$e^{-f_n}=E[e^{-f}|V_n] \;\;\text{and} \;\;e^{-g_n}=E[e^{-g}|V_n]$$ 
		where  $V_n$ is generated by $\{\delta{e_1},\dots,\delta{e_n}\}$ and $\{e_i,\ i\geq 1\}$ is an orthogonal basis of $H.$ Let $(\varphi_n,\psi_n)$ be the Monge potentials associated with Monge-Kantorovitch problem $(\rho_n,\nu_n),$ where $$\text{d}\rho_n=e^{-f_n}\text{d}\beta \text{ and } \text{d}\nu_n=e^{-g_n}\text{d}\beta$$
		and $\beta$ is the standard Gaussian measure on $\mathbb{R}^n$. Then, $(\varphi_n)$ converges to $\varphi$ in $\mathbb{D}_{2,1}(\rho)$, $(\psi_n)$ converges to $\psi$ in $L^1(\nu)$ and $(\nabla\psi_n)$ converges to $\nabla\psi$ in $L^2(\nu).$
	\end{lemma}
	\begin{proof}
		This lemma is proved in \cite[Thm.1.1]{fang}
	\end{proof}
	
As the next step, we will define  smooth functions $f_{m}$ and $g_{m}$ in the following lemma using the Ornstein-Uhlenbeck semigroup  to approximate functions $f$ and $g$ given in the finite dimension. The sequences $(f_{m})$ and $(g_{m})$ will be used later to approximate the sequences of  Lemma \ref{lemma1}, which are in finite dimension. 
\begin{lemma}\label{lemma2}
Let $\beta$ be the standard Gaussian measure on $\mathbb{R}^n$, $f\in\mathbb{D}_{2,1}(\beta)$ and $g\in\mathbb{D}_{2,1}(\beta)$ such that
\begin{align*}
	\int_{\mathbb{R}^d}|\nabla f|^2e^{-f}\;{d}\beta<\infty
	\text{ and }\int_{\mathbb{R}^d}|\nabla g|^2e^{-g}\;{d}\beta<\infty. 
\end{align*}
Let $(\varphi,\psi)$ be the Monge potentials associated with the Monge-Kantorovitch problem $\Gamma(\rho,\nu),$ where $$\text{d}\rho=e^{-f}\text{d}\beta \text{ and } \text{d}\nu=e^{-g}\text{d}\beta.$$ Define $f_m$ and $g_m$ as $$e^{-f_m}=Q_\frac{1}{m}e^{-f} \;\;\text{and} \;\;e^{-g_m}=Q_\frac{1}{m}e^{-g}$$ where $(Q_t, t\geq0)$ is the Ornstein-Uhlenbeck semigroup on $\mathbb{R}^d$. Let $(\varphi_m,\psi_m)$ be the Monge potentials associated with Monge-Kantorovitch problem $\Gamma(\rho_m,\nu_m),$ where $$\text{d}\rho_m=e^{-g_m}\text{d}\beta \text{ and }  \text{d}\nu_m=e^{-g_m}\text{d}\beta.$$ Then, $(\varphi_m)$ converges to $\varphi$ in $\mathbb{D}_{2,1}(\rho)$, $(Q_\frac{1}{m}\psi_m )$ converges to $\psi$ in $L^1(\nu)$ and $(Q_\frac{1}{m}\nabla\psi_m)$ converges to $\nabla\psi$ in $L^2(\nu).$
	\end{lemma}
	\begin{proof}
		Let $\gamma_m$ and $\gamma$ be solutions of Monge-Kantorovich problems for  $(\rho_m,\nu_m) $ and $(\rho,\nu)$, respectively. Then, we have
		\begin{align*}
		\int_{\mathbb{R}^d}|\nabla\varphi_m|^2e^{-f_m}\text{d}\beta&=d_2^2(\rho_m,\nu_m)
		\leq4\left(H(\rho_m\mid\beta)+ H(\nu_m\mid \beta)\right) \\
		&\leq4\left(H(\rho\mid\beta)+ H(\nu\mid \beta)\right) \\
		&\leq2\left(	\int_{\mathbb{R}^n}|\nabla f|^2e^{-f}\;{d}\beta+ 	\int_{\mathbb{R}^n}|\nabla g|^2e^{-g}\;{d}\beta \right)
		<\infty
		\end{align*}
		where we have used Talagrand's inequality in the first line, Jensen's inequality in the second and Logarithmic Sobolev  inequality in the third. Therefore, we have
		\begin{align*}
		\sup_{m}\int_{\mathbb{R}^d}|\nabla\varphi_m|^2e^{-f_m}\text{d}\beta<\infty
		\end{align*}
		Using Poincar{\'e} inequality, we get
		\begin{align*}
		\int_{\mathbb{R}^d}|\varphi_m-E_{\rho_m}[\varphi_m]|^2e^{-f_m}\text{d}\beta\leq\int_{\mathbb{R}^d}|\nabla\varphi_m|^2e^{-f_m}\text{d}\beta.
		\end{align*}
		Then, we replace $\varphi_m$ with $\varphi_m-E_{\rho_m}[\varphi_m]$ and use the fact that $\rho_m$ converges to $\rho$ weakly to get
		\begin{align*}
		\sup_{m}\|\varphi_m\|_{\mathbb{D}_{2,1}(\rho)}<\infty   \nnumber{bounded_sup}
		\end{align*}
		which implies $(\varphi_m)$ converges weakly in ${\mathbb{D}_{2,1}(\rho})$. 
As in the proof of \cite[Thm.4.1]{Feyel_Ustunel_2004},  we have
		\begin{align*}
		\varphi_m(x)+\psi_m(y)+\frac{1}{2}|x-y|^2\geq0.
		\end{align*} Applying Ornstein–Uhlenbeck semigroup with respect to $x$ and then $y$, we get 
		\begin{align*}
		Q_\frac{1}{m} \varphi_m(x)+Q_\frac{1}{m}\psi_m(y)+\frac{1}{2}Q_\frac{1}{m}\left( Q_\frac{1}{m}(|x-y|^2) \right) \geq0
		\end{align*}in view of the
positivity  improving property of Ornstein–Uhlenbeck semigroup. Observe that
		\begin{align*}
		&\lim\limits_{m}\int \left[Q_\frac{1}{m} \varphi_m(x)+Q_\frac{1}{m}\psi_m(y)+\frac{1}{2}Q_\frac{1}{m}Q_\frac{1}{m}(|x-y|^2|)\right]\;\text{d}\gamma \\&=\lim\limits_{m} \int Q_\frac{1}{m} \varphi_m(x)\;\text{d}\rho +\int Q_\frac{1}{m} \psi_m(y)\;\text{d}\nu+\int \frac{1}{2}Q_\frac{1}{m}Q_\frac{1}{m}(|x-y|^2|)\;\text{d}\gamma \nnumber{limit}
		\end{align*}
		We want to show that $Q_\frac{1}{m}Q_\frac{1}{m}(|x-y|^2|)$ converges to $|x-y|^2$ in $L^1(\gamma).$ 	Observe that
		\begin{align*}
		Q_\frac{1}{m}Q_\frac{1}{m}(|x-y|^2|)=\int\int\left|e^\frac{-1}{m}x+\sqrt{1-e^\frac{-2}{m}}t-e^\frac{-1}{m}y- \sqrt{1-e^\frac{-2}{m}}z\right|^2\; \text{d}\beta(z)\;\text{d}\beta(t).
		\end{align*} If we integrate with respect to $\gamma,$ we get 
		\begin{align*}
		\int Q_\frac{1}{m}Q_\frac{1}{m}(|x-y|^2|)\; \text{d}\gamma=\int\int \int &\left|e^\frac{-1}{m}x-e^\frac{-1}{m}y \right|^2 ({1-e^\frac{-2}{m}})|t|^2\\&+({1-e^\frac{-2}{m}})|z|^2\; \text{d}\beta(z)\;\text{d}\beta(t)\; \text{d}\gamma(x,y).
		\end{align*}
		It is easy to see that 
		\begin{align*}
		\lim\limits_{m\to\infty}\int({1-e^\frac{-2}{m}})|z|^2\; \text{d}\beta(z)=	\lim\limits_{m\to\infty}\int({1-e^\frac{-2}{m}})|t|^2\; \text{d}\beta(t)=0
		\end{align*}
		and 
		\begin{align*}
		\left|e^\frac{-1}{m}x+\sqrt{1-e^\frac{-2}{m}}t-e^\frac{-1}{m}y- \sqrt{1-e^\frac{-2}{m}}z\right|^2\leq C\left( |x|^2+|y|^2+|z|^2+|t|^2 \right).
		\end{align*}
		Moreover, from a version of Young's inequality given in \cite[pg. 15]{Beckenbach_1961}, we get
		\begin{align*}
		&\int |y|^2\;\text{d}\gamma(x,y)\leq \int e^{\alpha|y|^2}\text{d}\beta(y)+\frac{1}{\alpha}H(\nu|\beta)\\
		&\int |x|^2\;\text{d}\gamma(x,y)\leq \int e^{\alpha|x|^2}\text{d}\beta(x)+\frac{1}{\alpha}H(\rho|\beta).
		\end{align*}
		From dominated convergence theorem, $Q_\frac{1}{m}Q_\frac{1}{m}(|x-y|^2|)$ converge to $|x-y|^2$ in $L^1(\gamma).$ Therefore,  Equation \eqref{limit} is equal to
		\begin{align*}
		\lim\limits_{m} \int Q_\frac{1}{m} \varphi_m(x)\;\text{d}\rho &+\int Q_\frac{1}{m} \psi_m(y)\;\text{d}\nu+\int \frac{1}{2}|x-y|^2\ d\gamma 
		\\ &=\lim_{m}\int \varphi_m(x)\;\text{d}\rho_m +\int\psi_m(y)\;\text{d}\nu_m+\frac{1}{2}d_2^2(\rho,\nu) \\
		&=\lim_{m}-\frac{1}{2}d_2^2(\rho_m,\nu_m)+\frac{1}{2}d_2^2(\rho,\nu)\\
		&\leq \lim\limits_{m}-\frac{1}{2}d_2^2(\rho_m,\rho)+\frac{1}{2}d_2^2(\nu_m,\nu) \nnumber{distence_mimit}.
		\end{align*}
		 Observe that
		\begin{align*}
		\lim_m \int_{\mathbb{R}^d}|x|^2\;\text{d}\rho_m&=\lim_m \int_{\mathbb{R}^d}|x|^2 Q_\frac{1}{m}e^{-f}\;\text{d}\beta\\
		&=\lim_m \int_{\mathbb{R}^d}(Q_\frac{1}{m}\lambda|x|^2) \frac{e^{-f}}{\lambda}\;\text{d}\beta= \int_{\mathbb{R}^d}|x|^2 \;\text{d}\rho
		\end{align*}
		where we have used weak convergence of $\rho_m$ to $\rho$ in the first line and dominated convergence theorem in second. Indeed, Young's inequality  implies 
		$Q_\frac{1}{m}(\lambda|x|^2) \frac{e^{-g}}{\lambda}\leq \exp({Q_\frac{1}{m}\lambda|x|^2})+H(\rho|\beta)$ and for $\lambda<\frac{1}{2}$ Jensen's inequality yields
		\begin{align*}
		\int \exp({Q_\frac{1}{m}\lambda|x|^2})\text{d}\beta\leq \int Q_\frac{1}{m}\exp({\lambda|x|^2})
		= \int \exp({\lambda|x|^2})\text{d}\beta<\infty.
		\end{align*}
		Similarly,	$\lim \int_{\mathbb{R}^d}|y|^2\;\text{d}\nu_m= \int_{\mathbb{R}^d}|y|^2 \;\text{d}\nu.$ So $d_2(\rho_m,\rho)\to0$ and $d_2(\nu_m,\nu)\to0$ by using  \cite[Lem.8.3]{freedman}. 
		Combining this result with Equation \eqref{distence_mimit}, we get $(Q_\frac{1}{m} \varphi_m(x)+Q_\frac{1}{m}\psi_m(y)+\frac{1}{2}Q_\frac{1}{m}Q_\frac{1}{m}(|x-y|^2) )$ converges to 0 in $L^1(\gamma)$ and it is uniformly integrable with respect to $\gamma$. Moreover, $(Q_\frac{1}{m}\varphi_m)$ is uniformly integrable, so $(Q_\frac{1}{m}\psi_m)$ is uniformly integrable. Let $a'$ and  $b'$ be weak limit points of $(Q_\frac{1}{m}\varphi_m)$ and $(Q_\frac{1}{m}\psi_m).$ The Cesaro means
		\begin{align*}
		Q_\frac{1}{m}\varphi_m'=\frac{1}{n}\sum_{i=1}^{n}Q_\frac{1}{m}\varphi_m\;\;\text{ and }\;\;Q_\frac{1}{m}\psi_m'=\frac{1}{n}\sum_{i=1}^{n}Q_\frac{1}{m}\psi_m.
		\end{align*} 
	 converge to $a'$, up to a subsequence, and $b'$  in $L^1(\gamma)$, respectively. A a result, we have
		\begin{align*}
		a'+b'+\frac{1}{2}|x-y|^2= 0\quad\gamma\text{-a.s.}
		\end{align*}
		Let $a(x)=\limsup\varphi_m'(x)$ and $b(y)=\limsup\psi_m'(y).$ Then $a'=a$ and $b'=b$ $\gamma$-a.s. and  
		\begin{align*}
		a(x)+b(y)+\frac{1}{2}|x-y|^2&\geq0,
		\quad\text{for any $(x,y)\in\mathbb{R}^d\times\mathbb{R}^d$}\\
		a(x)+b(y)+\frac{1}{2}|x-y|^2&=0,\quad
		\text{$\gamma-$a.s.}
		\end{align*}
		By uniqueness of solutions, we get $a=\varphi$ and $b=\psi.$ Therefore, we deduce that $(Q_\frac{1}{m}\varphi_m)$ converges weakly to $\varphi$ in $L^1(\gamma)$.
		
		On the other hand, from \eqref{bounded_sup}, there exists $\varphi'\in L^2(\beta)$ such that $(\varphi_m)$ converges weakly to $\varphi'$ in $L^2(\rho)$ and in $L^2(\gamma).$ For $h\in L^2(\gamma),$ we have
		
		$$ \int\left(\varphi_m-\varphi\right)h\;d\gamma=\int\left(\varphi_m-Q_\frac{1}{m}\varphi_m\right)h\;d\gamma +\int\left(Q_\frac{1}{m}\varphi_m-\varphi\right)h\;d\gamma.$$
		Since both integrals on the right hand side converge to zero as $m\to\infty$, $(\varphi_m)$ converges weakly to $\varphi$ in $L^2(\gamma)$ and in $L^2(\rho)$.  Moreover,  $\lim \E_\rho [|\nabla \varphi_m|^2]=\E_\rho[|\nabla \varphi|^2] $, which implies $(\varphi_m)$ converges to $\varphi$ in $\mathbb{D}_{2,1}(\rho).$
		Similarly, $(Q_\frac{1}{m}\psi_m)$ converges to $\psi$ in $L^1(\nu)$, and since $\nabla$ is closable in $L^p(\nu)$ for $p\geq 1$, $(\nabla Q_\frac{1}{m}\psi_m)$  converges weakly to $\nabla\psi$ in $L^2(\nu).$ In addition, we have 
		\begin{align*}
		\lim\limits_{m}\E_{\nu_m}\left[ |\nabla\psi_m|^2\right] &=\lim\limits_{m}d_2^2(\rho_m,\nu_m)=d_2^2(\rho,\nu)=\E_\nu\left[|\nabla\psi|^2 \right] 
		\end{align*}
		which implies that $Q_\frac{1}{m}\nabla\psi_m$ converges to $\nabla\psi$ in $L^2(\nu)$.
	\end{proof}
	Now, let $\rho$ and $\nu$ be probability measures on $\R^n$ defined by 
	\begin{align*}\nnumber{abs_measures}
	d\rho=F\;d\beta ,    \quad \quad 	d\nu=G\;d\beta.
	\end{align*}
	where  $\beta$ is a Gaussian measure on $\R^n$ and $F, G\in L^1(\beta).$ Suppose that $\nu$ is also absolutely continuous with respect to $\rho$ with $\frac{d\nu}{d\rho}=L.$ 
	Define 
	\begin{align*}
	F_k=\frac{\theta_k F}{\E[\theta_k F]},  \quad \quad  G_k=\frac{\theta_k G }{\E[\theta_k G]}
	\end{align*}
	where $\theta_k\in C_c^\infty(\R^n)$ is a smooth function with compact support satisfying  $\theta_k(x)=0$ if $|x|\geq k,$ $\theta_k(x)=1$ if $|x|\leq k-1$ with $0\leq \theta_k\leq1$ for each $k\in\N$, and $\sup_k|\frac{(\nabla\theta_k)^2}{\theta_k}|\leq 1$. 
	Consider probability measures $\rho_k$ and $\nu_k$ given by
	\begin{align*}\nnumber{compact_abs_measures}
	d\rho_k=F_k\;d\beta,  \quad \quad  d\nu_k=G_k\;d\beta.
	\end{align*}
	It is easy to see  that $\nu_k$ is absolutely continuous with respect to $\rho_k$ with Radon-Nikodym derivative ${d\nu_k}/{d\rho_k}=L_k, $ where 
	$$L_k= \begin{cases} 
	\frac{G_k}{F_k} & \text{on }\{\theta_k\neq0 \} \\
	0 & \text{on }\{\theta_k=0 \}
	\end{cases}. $$ 
	Observe that $L_k={b_{k}}/{a_{k}}L$ on the set $\{\theta_k\neq0 \},$ where $a_k={1}/{\E[\theta_k F]}$ and $b_k={1}/{\E[\theta_k G]}$.  Therefore, the relative entropy of $\nu_k$ with respect to  $\rho_k$, 
	$H(\nu_k|\rho_k)$, is finite if and only if $H(\nu|\rho)$ is finite.
	Note that, $F_{k}$ and $G_{k}$ are bounded, being  continuous functions with compact support. Using the following lemma, we will be able to approximate the smooth density functions of Lemma \ref{lemma2} by smooth and bounded sequences.
	
	\begin{lemma}\label{lemma3}
		Let $(\rho,\nu)$ and $(\rho_k,\nu_k)$  be probability measures on $\R^n$ defined as in \eqref{abs_measures} and \eqref{compact_abs_measures}, respectively, with $H(\nu|\rho)<\infty.$ Let $(\varphi,\psi)$ and $(\varphi_k,\psi_k)$ be the Monge potentials associated with Monge-Kantorovitch problems for $(\rho,\nu)$ and $(\rho_k,\nu_k),$ respectively, with quadratic cost. Then $(\varphi_k)$ converges to $\varphi$ in $\mathbb{D}_{2,1}(\rho)$, $(\theta_k\psi_k)$ converges to $\psi$ in $L^1(\nu)$ and $(\sqrt{\theta_k}\nabla\psi_k)$ converges to $\nabla\psi$ in $L^2(\nu).$
	\end{lemma}
	\begin{proof}
		Let $\gamma$ and $\gamma_k$ be the optimal transport plans to the Monge-Kantorovich problems for $(\rho,\nu)$ and $(\rho_k,\nu_k).$ From \cite[Prop.3.1]{Brenier}, replacing $\varphi_k$ with $\varphi_k-E[\varphi_k]$ and $\psi_k$ with $\psi_k+E[\varphi_k]$, we have
		\begin{align*}
		F_k(x,y)=\varphi_n(x)+\psi_n(y)+\frac{1}{2}|x-y|^2\geq0\nnumber{com_F_n}
		\end{align*}
		for all $x,y\in\R^n$.
		Since $0\leq \theta_k$, we have 
		\begin{align*}
		0&\leq\int_{\R^n\times \R^n} \theta_k(y) F_k(x,y) \ d\gamma(x,y)\\
		&=\int \theta_k(y)\varphi_k(x)\ d\gamma+\int \theta_k(y)\psi_k(y)\ d\gamma +\frac{1}{2}\int\theta_k(y)|x-y|^2\ d\gamma. \nnumber{three_parts}
		\end{align*}	
	In the second integral above, we have
		\begin{align*}
		\int \theta_k(y)\psi_k(y)\ d\gamma
		&=\int \theta_k(y)\psi_k(y)\ d\nu\\
		&=b_k^{-1}\int \psi_k(y)\;d\nu_k(y)\\
		&=b_k^{-1}\int \psi_k(y)\;d\gamma_k(y)\\
		&=b_k^{-1}\int -\varphi_k(x)-\frac{1}{2}|x-y|^2\;d\gamma_k(x,y)\\
		&=b_k^{-1}\left[ \int -\varphi_k(x)\;\rho_k(x) -\frac{1}{2}\int|x-y|^2\;d\gamma_k(x,y)\right] \\
		&=\frac{a_k}{b_k} \int -\theta_k(x)\varphi(x)\;d\rho-\frac{1}{2b_k}\int |x-y|^2\;d\gamma_k(x,y) \nnumber{second}
		\end{align*}	If we substitute \eqref{second} into \eqref{three_parts}, we get
		\begin{align*}
		\int\theta_k(y) F_k(x,y)\ d\gamma=&\int \left( \theta_k(y)-\frac{a_k}{b_k}\theta_k(x)\right) \varphi_k(x)\ d\gamma \\
		&+\frac{1}{2}\int \theta_k(y)|x-y|^2\;d\gamma(x,y) -\frac{1}{2b_k}\int |x-y|^2\;d\gamma_k(x,y)\\
		&=:\Romannum{1}_k+\Romannum{2}_k-\Romannum{3}_k 
		\end{align*}
		Sequences $(a_k )$ and $(b_k)$  given by $a_k=\E[\theta_k F]^{-1}$ and  $b_k=\E[\theta_k G]^{-1}$ converge to 1 decreasingly. Therefore, for given $\epsilon>0$ there exists $0<k_\epsilon$ such that $0<\frac{b_k}{a_k}\leq \frac{1}{1-\epsilon}$ for each $k\geq k_\epsilon$. Hence, we have
		\begin{align*}
		\lim_{k} \int |\varphi_k|^2 \;d\gamma&=\lim_{k} \int |\varphi_k|^2 \;d\beta
		 \leq\lim_{k} \int |\nabla \varphi_k|^2\;d\rho \\
		& \leq\lim_{k} 2 H(\nu_k| \rho_k)
		\leq 2\int \left( \frac{b_k}{a_k}L \log L + \frac{b_k}{a_k}L \log \frac{b_k}{a_k} \right)\;d\rho\\ \nnumber{varhi_bound}
		&\leq \frac{2}{1-\epsilon}H(\nu|\rho)+\frac{2}{1-\epsilon}<\infty
		\end{align*}
		for  sufficiently small $\epsilon.$ We have used Poincar\'{e} inequality in the first line and Talagrand's inequality in the second. Therefore, from Cauchy–Schwarz inequality and dominated convergence theorem, $\Romannum{1}_k\to 0$ as $k\to \infty$.
		From monotone convergence theorem, $\Romannum{2}_k$ converges to $\frac{1}{2}\int |x-y|^2\;d \gamma.$
		Moreover,  $\Romannum{3}_k$ also converges to $\frac{1}{2}\int |x-y|^2\;d \gamma$. Indeed, Young's inequality implies that
		\begin{align*}
		&|x|^2\theta_k(x)F\leq\frac{1}{1-\epsilon}\left( e^{\epsilon|x|^2}+\frac{1}{\epsilon}F\log F\right) \\
		&|y|^2\theta_k(y)G\leq\frac{1}{1-\epsilon}\left( e^{\epsilon|y|^2}+\frac{1}{\epsilon}G\log G\right) .
		\end{align*}
		In other words, the sequences $(|x|^2\theta_k(x)F,\;k\geq k_\epsilon)$ $(|y|^2\theta_k(y)G,\;k\geq k_\epsilon)$ are bounded by a $\beta-$integrable function. Using dominated convergence theorem, we get
		\begin{align*}
		\lim\limits_{k\to\infty} \int |x|^2F_k\;d\beta=\lim\limits_{k\to\infty} \int |x|^2\;d\rho_k=\int |x|^2\;d\rho\\
		\lim\limits_{k\to\infty} \int |y|^2G_k\;d\beta=\lim\limits_{k\to\infty} \int |y|^2\;d\nu_k=\int |y|^2\;d\nu.
		\end{align*}
		Combining the above with \cite[Lem.8.3]{freedman} yields
		\begin{align*}
		\lim\limits_{k\to\infty}\frac{1}{2b_k}\int |x-y|^2\;d\gamma_k(x,y)&=\lim\limits_{k\to\infty} \frac{1}{2b_k}d_2^2(\rho_k,\nu_k)\\
		&=\frac{1}{2}d_2^2(\rho,\nu)
		=\frac{1}{2}\int |x-y|^2\;d\gamma(x,y).
		\end{align*}
Therefore, we get
		\begin{align*} 
		\lim\limits_{k\to\infty} 	\int \theta_k(y) F_k(x,y)\ d\gamma =0. \nnumber{sifir}
		\end{align*}
		The sequence $(\varphi_k)$ is bounded in $\mathbb{D}_{2,1}(\rho)$. Therefore, it converges weakly to some $\varphi'\in \mathbb{D}_{2,1}(\rho)$ up to a subsequence. Moreover, the sequence $(\theta_k(y)\varphi_k(x))$
		is uniformly integrable and for any $h\in L^\infty(\gamma)$
		\begin{align*}
		\int \left( \theta_k(y)\varphi_k(x) - \varphi' (x)\right) h\;d\gamma= \int \left(\theta_k(y) - 1 \right) \varphi_k(x) h\;d\gamma+\int \left( \varphi_k(x) - \varphi' (y)\right) h\;d\gamma.
		\end{align*}
		The second integral on the right hand side converges to zero since $(\varphi_k)$ converges to $\varphi'$ weakly and the first integral converges to zero since $(\theta_k-1)$ converges to zero in $L^2(\gamma)$. 
Therefore,  $(\theta_k\varphi_k)$ converges weakly to $\varphi'$ in $L^2(\gamma)$.
 On the other hand,   $(\theta_k F_k)$ converges to zero in $L^1(\gamma)$ as given in \eqref{sifir},  so it is uniformly integrable. Since $(\theta_k\varphi_k)$ is uniformly integrable,   $(\theta_k(y)\psi_k(x))$ is also uniformly integrable and  converge weakly to some $\psi'$ in $L^1(\gamma)$ up to a subsequence. 
There exists a further subsequence such that the Cesaro means  
		\begin{align*}
		\varphi_k'=\frac{1}{k}\sum_{i=1}^{k}\theta_i(y)\varphi_i(x)\;\;\text{ and }\;\;\psi_n'=\frac{1}{k}\sum_{i=1}^{k}\theta_i\psi_i.
		\end{align*} 
converge to $\varphi'$ in $L^2(\gamma)$ and $\psi'$  in $L^1(\gamma)$, respectively. Let $a(x)=\limsup\varphi_k'(x)$ and $b(y)=\limsup\psi_k'(y).$ Then, $\varphi'=a$ and $\psi'=b$ $\gamma$-a.s.  
		Therefore, $\varphi'=a$ and $\psi'=b$ $\gamma$-a.s. and  
		\begin{align*}
		a(x)+b(y)+\frac{1}{2}|x-y|^2&\geq0,
		\quad\text{for any $(x,y)\in\mathbb{R}^d\times\mathbb{R}^d$}\\
		a(x)+b(y)+\frac{1}{2}|x-y|^2&=0,\quad
		\text{$\gamma$-a.s.}
		\end{align*}
		By uniqueness of solutions, we have $\varphi'=\varphi$ and $\psi'=\psi.$  Combining this result with $\lim \E_\rho[|\nabla\varphi_k|^2]=\E_\rho[|\nabla\varphi|^2]$, we see that $(\varphi_k)$ converges strongly to $\varphi$ in $\mathbb{D}_{2,1}(\rho)$.   Finally, since $(\theta_k F_k)$ converges to zero and  $(\frac{1}{2}\theta_k |x-y|^2)$ converges to $\frac{1}{2}|x-y|^2$ in $L^1(\gamma)$,  $\theta_k\psi_k$ converges to $\psi$  in $L^1(\gamma)$.
		The only remaining part is to show the convergence of the sequence $(\sqrt{\theta_k} \nabla\psi_k)$. Observe that
		\begin{align*}
		\int \theta_k|\nabla \psi_k|^2\ d\gamma&=\int \theta_k|\nabla\psi_k|^2\ d\nu
		=\int b_k^{-1}|\nabla\psi_k|^2 \;d\nu_k(y)\\
		&=\int  b_k^{-1}|\nabla\varphi_k|^2\;d\rho_k(x)
		=\int \frac{a_k}{b_k}\theta_k|\nabla\varphi_k|^2\ d\rho(x)
		\end{align*}
		If we take the limit of both sides, we get
		\begin{align*}
		\lim\limits_{k\to \infty}\int \theta_k|\nabla \psi_k|^2\ d\gamma
		&=\int |\nabla\varphi|^2\;d\rho(x)=\int |\nabla\psi|^2\;d\nu(y)\nnumber{psiconverge}.
		\end{align*} 
		Next, we will show that the sequence $(\sqrt{\theta_k}\nabla\psi_k)$ converges weakly to $\nabla\psi.$ Let $\xi$ be a bounded smooth vector field, then
		\begin{align*}
		\int \langle\sqrt{\theta_k}\nabla\psi_k,\xi\rangle\;d\gamma=\int_{\{\theta_k<c\}}\langle\sqrt{\theta_k}\nabla\psi_k,\xi\rangle\;d\gamma+\int_{\{\theta_k\geq c\}}\langle\sqrt{\theta_k}\nabla\psi_k,\xi\rangle\;d\gamma.\nnumber{twoparts}
		\end{align*}
Equation \eqref{psiconverge} implies that the sequence $(\sqrt{\theta_k}\nabla\psi_k)$ is uniformly integrable with respect to $\gamma$. Combining this with the boundedness of $\xi$ and the fact that $\theta_k$ converges to $1$ imply that the first integral converges to $0$ for $c<1.$ On the other hand, the second integral in $\eqref{twoparts}$ can be written as 
		\begin{align*}
		\int_{\{\theta_k\geq c\}}\langle\sqrt{\theta_k}\nabla\psi_k,\xi\rangle\;d\nu&=b^{-1}_{ k}\int_{\{\theta_k\geq c\}} \frac{\langle\nabla\psi_k,\xi\rangle}{\sqrt{\theta_k}}\;d\nu_k\\
		&=b^{-1}_{ k}\int_{\{\theta_k\circ T_k\geq c\}} -\frac{\langle\nabla\varphi_k,\xi\circ T_k\rangle}{\sqrt{\theta_k\circ T_k}}\;d\rho_k\\
		&=\frac{a_k}{b_k}\int_{\{\theta_k\circ T_k\geq c\}} -\frac{\langle\nabla\varphi_k,\xi\circ T_k\rangle}{\sqrt{\theta_k\circ T_k}}\theta_k\;d\rho
		\end{align*}
where $T_k=I_{\R^n}+\nabla\varphi_k$ is the optimal transport map of the Monge problem.
	Since the sequence $(L_k)$ is uniformly  integrable, the sequence $(T_k)$ is equi-concentrated on a compact set. Hence, $\lim\limits_{k\to \infty}\theta_k\circ T_k=1$ and $\lim\limits_{k\to \infty}\xi\circ T_k=\xi\circ T$ in $\rho$-probability. By using dominated convergence theorem and $L^2(\rho)$-boundedness of the sequence of $(\nabla\varphi_k)$, we get
		\begin{align*}
		\lim\limits_{k\to \infty}\frac{a_k}{b_k}\int_{\{\theta_k\circ T\geq c\}} -\frac{\langle\nabla\varphi_k,\xi\circ T_k\rangle}{\sqrt{\theta_k\circ T}}\theta_k\;d\rho
		&=\int- \langle\nabla\varphi,\xi\circ T\rangle\;d\rho\\
		&=\int (\nabla\psi_k\circ T,\xi\circ T)\;d\rho\\
		&=\int (\nabla\psi,\xi)\;d\gamma.
		\end{align*} 
		As a result, the limit of  \eqref{twoparts} is 
		\begin{align*}
		\lim\limits_{k\to \infty} \int \langle\sqrt{\theta_k}\nabla\psi_k,\xi\rangle\;d\gamma=\int \langle\nabla\psi_k,\xi\rangle\;d\gamma
		\end{align*}
		which implies that the sequence $\sqrt{\theta_k}\nabla\psi_k$ converges weakly to $\nabla\psi$ with respect to $\gamma.$ In view of \eqref{psiconverge},   $(\sqrt{\theta_k}\nabla\psi_k)$ converges to $\nabla\psi$ in  $L^2(\nu)$. 
	\end{proof}

	In the next lemma, we will  use strictly positive sequences to approximate  density functions such as those used in Lemma \ref{lemma3}. 
	
	\begin{lemma}\label{lemma4}
		Let $(\rho,\nu)$  be probability measures on $\R^n$ defined as in \eqref{abs_measures}  with $H(\rho\mid \beta)<\infty$. Define for each $\epsilon>0$
		\begin{align*}
		d\rho_\epsilon=\frac{F+\epsilon}{1+\epsilon}\;d\beta\quad\text{and}\quad
		d\nu_\epsilon=\frac{G+\epsilon}{1+\epsilon}\;d\beta.
		\end{align*}
		Let $(\varphi,\psi)$ and $(\varphi_\epsilon,\psi_\epsilon)$ be the Monge potentials associated with Monge-Kantorovitch problems $(\rho,\nu)$ and $(\rho_\epsilon,\nu_\epsilon)$, respectively with quadratic cost. Then, as $\epsilon$ goes to zero $(\varphi_\epsilon,\;\epsilon>0)$ converges to $\varphi$ in $\mathbb{D}_{2,1}(\rho)$, $(\psi_\epsilon,\;\epsilon>0)$ converges to $\psi$ in $L^1(\nu)$ and $(\nabla\psi_\epsilon,\;\epsilon>0)$ converges to $\nabla\psi$ in $L^2(\nu).$
	\end{lemma}
	\begin{proof}
		Let $\gamma$ and $\gamma_\epsilon$ be the optimal transport plans of Monge-Kantorovich problems for $(\rho,\nu)$ and $(\rho_\epsilon,\nu_\epsilon)$, respectively. If we  replace $\varphi_\epsilon$ with $\varphi_\epsilon-\E_{\rho_\epsilon}[\varphi_\epsilon]$ and $\psi_\epsilon$ with $\psi_\epsilon+\E_{\rho_\epsilon}[\varphi_\epsilon]$, we have
		\begin{align*}
		F_\epsilon(x,y)&=\varphi_\epsilon(x)+\psi_\epsilon(y)+\frac{1}{2}|x-y|^2\geq0 \quad\text{for any } x,y\in\R^n,\\
		F_\epsilon(x,y)&=\varphi_\epsilon(x)+\psi_\epsilon(y)+\frac{1}{2}|x-y|^2=0 \quad\gamma_\epsilon \text{-a.s.} \nnumber{epsilon_F}
		\end{align*}
		for all $x,y\in\R^n$.
		Observe that 
		\begin{align*}
		0&\leq\int F_\epsilon(x,y)\ d\gamma+\epsilon\int  \varphi_\epsilon(x)\ d\beta(x)+\epsilon\int  \psi_\epsilon(y)\ d\beta(y)\\
		&=(1+\epsilon)\int  \varphi_\epsilon(x)\ d\rho_\epsilon(x)+(1+\epsilon)\int  \psi_k(y)\ d\nu_\epsilon(y) +\frac{1}{2}\int|x-y|^2\ d\gamma\\
		&=(1+\epsilon)\int  \psi_k(y)\;\nu_\epsilon(y) +\frac{1}{2}\int|x-y|^2\ d\gamma  \nnumber{epsilon_three_parts}.
		\end{align*}
		Note that last inequality follows from the fact that $\E_{\rho_\epsilon}[\varphi_\epsilon]=0$. Using Young's inequality and \cite[Lem.8.3]{freedman}, we get
		\begin{align*}
		\lim\limits_{\epsilon\to 0} (1+\epsilon)\int  \psi_\epsilon\ d \nu_\epsilon&=\lim\limits_{\epsilon\to 0} \frac{1+\epsilon}{2}\int|x-y|^2\;d\gamma_\epsilon\\
		&=\lim\limits_{\epsilon\to 0}-d_2^2(\rho_\epsilon,\nu_\epsilon)
		=-d_2^2(\rho,\nu)
		=-\frac{1}{2}\int|x-y|^2\;d\gamma
		\end{align*}
Therefore, we have
		\begin{align*}
		\lim\limits_{\epsilon\to 0} \int F_\epsilon(x,y)\ d\gamma+\epsilon\int  \varphi_\epsilon(x)\;d\beta(x)+\epsilon\int  \psi_\epsilon(y)\;d\beta(y)=0.
		\end{align*}
		On the other hand, $F_\epsilon(x,x)=\varphi_\epsilon(x)+\psi_\epsilon(x)\geq0$. Hence, $\epsilon\int  \varphi_\epsilon\;\beta+\epsilon\int  \psi_\epsilon\;\beta\geq0$, which implies
		\begin{align*}
		\lim\limits_{\epsilon\to 0} \int F_\epsilon(x,y)\ d\gamma=0
		\end{align*}
		and we conclude that $(F_\epsilon,\;\epsilon>0)$ is uniformly integrable with respect to $\gamma.$ 
		Observe that 
		\begin{align*}
		\int|\nabla\varphi_\epsilon|^2\;\text{d}\rho_\epsilon&=d_2^2(\rho_\epsilon,\nu_\epsilon)\\
		&\leq4\left(H(\rho_\epsilon\mid\beta)+ H(\nu_\epsilon\mid\beta)\right).
		\end{align*}		
		Moreover,  $V_\epsilon=\frac{V}{1+\epsilon}+\frac{\epsilon}{1+\epsilon}$ and the convexity of the function  $x\to x\log x$ implies
		\begin{align*}
		H(\rho_\epsilon\mid \beta)=\int F_\epsilon\log F_\epsilon\;d\beta	\leq \frac{1}{1+\epsilon}H(\rho\mid\beta)
		\end{align*}
	Similarly,  $H(\nu_\epsilon|\beta)\leq\frac{1}{1+\epsilon}H(\nu|\beta)$ and
		\begin{align*}
		\sup_{\epsilon>0}\int|\nabla\varphi_\epsilon|^2\text{d}\rho_\epsilon<\infty.
		\end{align*}
Combining the above result with $\int|\nabla\varphi_\epsilon|^2\;\text{d}\rho_\epsilon\geq\frac{1}{1+\epsilon}\int|\nabla\varphi_\epsilon|^2\;\text{d}\rho$, we get 
$\sup_{\epsilon>0}\int|\nabla\varphi_\epsilon|^2\text{d}\rho<\infty$. Applying  Poincar\'{e} inequality yields
		$\sup_{\epsilon>0} \|\varphi_\epsilon\|_{\mathbb{D}_{2,1}(\rho)}<\infty$. Since the sequence $(\varphi_\epsilon, \epsilon>0)$  is bounded in $\mathbb{D}_{2,1}(\rho)$, it is also uniformly integrable with respect to $\rho$ and $\gamma.$ We also know that the sequence $(F_\epsilon, \epsilon>0)$ is uniformly integrable. Hence, $(\psi_\epsilon,\;\epsilon>0)$ is uniformly integrable with respect to $\gamma.$  The rest of the proof goes along the proof of the previous lemmas. 
	\end{proof}
	
	\subsection{Regularity of Backward Monge Potential}\label{Regularity of Backward Monge Potential}  \label{s4.2}
	
	Our approach for proving regularity of backward Monge potential in a more general setting  will be to approximate the functions $f$ and $g$ of \eqref{conds} by  appropriate sequences  that will enable the use of  Proposition \ref{regulary_smooth_bakcward}. Explicitly, we first consider $f_n$ on $\mathbb{R}^n$ by $e^{-f_n }=\E[e^{-f} |V_n]$, where  $V_n$ is generated by $\{\delta{e_1},\dots,\delta{e_n}\}$ for an orthogonal basis $\{e_i,\ i\geq 1\}$  of $H$.  In the second step, $f_{nm}$ are chosen as smooth functions using the Ornstein-Uhlenbeck semigroup $(P_{1/m})$ by $e^{-f_{nm} }=P_{1/m} (e^{-f_n })$, for each $n$. Then, $F_{nmk}$  are continuous and compact supported functions given by $F_{nmk}=\frac{\theta_k e^{-f_{nm }}}{\E[\theta_k e^{-f_{nm }}]}$ for fixed $m$ and $k$. Finally,  a  strictly positive sequence given by 
	$e^{-f_{nmkl} }=\frac{F_{nmk}+\frac{1}{l}}{1+\frac{1}{l}}$
	is formed to be compatible with the density $e^{-f}>0$. Similarly, a sequence $(g_{nmkl}) $ is also defined.  As a result, forward Monge potential $\varphi_{nmkl}$ of the Monge-Kantorovich problem with respect to $(\rho_{nmkl},\nu_{nmkl})$ is smooth and Proposition \eqref{regulary_smooth_bakcward} is applicable. 
	We will work with a subsequence of $(\varphi_{nmkl})$, which has the form $(\varphi_{n,k_n,l_{k_n}, m_{l_{k_n}}})$ and can be extracted by applying the diagonal method three times. We will denote the corresponding  sequences with $(\varphi_n)$, $(f_n)$ and $(\rho_n)$.
	Similarly,  $(\psi_n)$ will be formed by the diagonal method and its indices can be matched with those of $(\varphi_n)$. We will work with this subsequence by relabeling as $(\varphi_n,\psi_n)$.  After this simplification, $(\varphi_n)$ converges to $\varphi$ in $\D_{2,1}(\rho),$ $(\theta_nP_{1/n}\psi_n)$ converges to $\psi$ in $L^1(\nu)$ and $(\sqrt{\theta_{{n}}} P_{1/n}\nabla\psi_n)$ converges to $\nabla\psi$ in $L^2(\nu)$.
	\begin{lemma}\label{lemma5} Let $f,g\in\mathbb{D}_{2,1}$ satisfy \eqref{conds} with $e^{-f} $ satisfying Poincar{\'e} inequality \eqref{poincare}. Then, we have
		\begin{align*} 
		\lim\limits_{n}  \nabla f_{n}=\nabla f,  \quad \quad 
		\lim\limits_{n}  \nabla g_{n}\circ T_{n}=\nabla g\circ T
		\end{align*}
		in $L^2(\rho,H)$.
	\end{lemma}
	
	\begin{proof} It will be convenient to use the multi-index sequence defined above as  $$ e^{-f_{nmkl}}=\frac{\theta_k e^{-f_{nm}}+\frac{1}{l}}{1+\frac{1}{l}}\; .$$
		We have
		\begin{align*}
		\E\left[ |\nabla f_{nmkl}|^2 e^{-f_{nmkl}} \right]&=4\E\left[ \left| \nabla e^{-f_{nmkl}/2} \right|^2  \right]\; = \;	\E\left[ \left| \frac{\nabla e^{-f_{nmkl}}}{e^{-f_{nmkl/2}}} \right|^2  \right]
		\\
		&= 	\E\left[1_{\{\theta_k=0\}}  \left| \frac{\nabla e^{-f_{nmkl}}}{e^{-f_{nmkl/2}}} \right|^2  \right]+\E\left[1_{\{0<\theta_k<1\}}  \left| \frac{\nabla e^{-f_{nmkl}}}{e^{-f_{nmkl/2}}} \right|^2  \right]\\ &\qquad\qquad+\E\left[1_{\{\theta_k=1\}}  \left| \frac{\nabla e^{-f_{nmkl}}}{e^{-f_{nmkl/2}}} \right|^2  \right]
		\\
		&=\Romannum{1}_{nmkl}+\Romannum{2}_{nmkl}+\Romannum{3}_{nmkl}
		\end{align*}
		On the set ${\{\theta_k=1\}}$, we have
		\begin{align*}
		\Romannum{3}_{nmkl}&=\frac{l}{l+1}\E\left[ 1_{\{\theta_k=1\}} \left| \frac{\nabla e^{-f_{nm}}}{\left( e^{-f_{nm}}+\frac{1}{l}\right)^\frac{1}{2} } \right|^2  \right]
		\\
		&\leq  \E\left[ 1_{\{\theta_k=1\}} \left| \nabla P_{\frac{1}{m}}(e^{-f_{n}})\right|^2  \left( P_{\frac{1}{m}} (e^{-f_{n}})\right)^{-1}  \right]
		\\
		&\leq  {e^{-2/m}}\E\left[ 1_{\{\theta_k=1\}} P_{\frac{1}{m}} \left(  \left|  \nabla e^{-f_{n}}\right|^2\right)   \left( P_{\frac{1}{m}} (e^{-f_{n}})\right)^{-1}  \right]
		\\
		&=  \E\left[ 1_{\{\theta_k=1\}}P_{\frac{1}{m}} \left( \left|  \nabla f_n  e^{-f_{n}}\right|^2\right)   \left( P_{\frac{1}{m}} (e^{-f_{n}})\right)^{-1}  \right]
		\\
		&=  \E\left[ 1_{\{\theta_k=1\}}P_{\frac{1}{m}} \left( \left|  \nabla f_n\right|^2  e^{-f_{n}}\right)  \frac{P_{\frac{1}{m}} (e^{-f_{n}})}{P_{\frac{1}{m}} (e^{-f_{n}})}  \right] \\
		&\leq         \E\left[ 1_{\{\theta_k=1\}} \frac{\left| \nabla e^{-f_n} \right|^2 }{e^{-f_n}} \right]\\
		&\leq         \E\left[ 1_{\{\theta_k=1\}} \frac{\left| \nabla \E[e^{-f}|V_n] \right|^2 }{\E[e^{-f}|V_n]} \right]
		\leq         \E\left[ 1_{\{\theta_k=1\}} \frac{\E[\left| \nabla e^{-f}\right|^2|V_n ] }{E[e^{-f}|V_n]} \right]\\
		&\leq         \E\left[ 1_{\{\theta_k=1\}} \E[\left| \nabla f\right|^2 e^{-f}|V_n ] \frac{\E[e^{-f}|V_n]}{E[e^{-f}|V_n]} \right]
		\leq         \E\left[  \left| \nabla f\right|^2 e^{-f} \right]
		\end{align*}
		On the set ${\{0<\theta_k<1\}}  $, we get
		\begin{align*}
		\Romannum{2}_{nmkl}&=\frac{l}{1+l}\E\left[1_{\{0<\theta_k<1\}}  \frac{\left( \nabla \theta_k e^{-f_{nm}}+\theta_k\nabla e^{-f_{nm}}\right)^2} {\theta_k e^{-f_{nm}}+\frac{1}{m}} \right]\\
		&\leq 2\E \left[ 1_{\{0<\theta_k<1\}}  \frac{|\nabla \theta_k|^2}{\theta_k} e^{-f_{nm}}\right]+ 2\E \left[ \frac{|\nabla e^{-f_{nm}}|^2}{e^{-f_{nm}}}\right]\\
		&\leq 2 \E \left[ 1_{\{0<\theta_k<1\}}  \frac{|\nabla \theta_k|^2}{\theta_k} e^{-f} \right]+2 \E \left[ \frac{|\nabla e^{-f}|^2}{e^{-f}}\right]
		\end{align*}
		Therefore, $	\Romannum{2}_{nmkl}$ is bounded uniformly. Moreover, $\lim\limits_{k}\mu({\{0<\theta_k<1\}})=0,$ which implies  $\lim \Romannum{2}_{nmkl}=0. $ When ${\{\theta_k=0\}}$, as the function $e^{-f_{nmkl}}$
		is constant, its derivative is zero and 
		$\lim\Romannum{1}_{nmkl}=0$ as well. 
		As a result, $(\nabla e^{-f_{nmkl}/2} )$ is uniformly bounded in $L^2(\mu,H)$ and 		\begin{align*}
	\limsup \E\left[ \left| \nabla e^{-f_{nmkl}/2} \right|^2  \right] 
	 \leq \E\left[ \left| \nabla e^{-f/2} \right|^2  \right].
\end{align*}
		On the other hand, $(f_{nmkl})$ converges to $f$ in $L^0(\mu)$, that is, in probability, $\E\left[ \left|  e^{-f_{nmkl}/2} \right|^2  \right]$ converges to $\E\left[ \left|  e^{-f/2} \right|^2  \right]$, which imply that $(e^{-f_n/2})$ converges to $e^{-f/2}$ in $L^2(\mu)$ and $(\nabla e^{-f_{nmkl}/2})$ converges weakly to $\nabla e^{-f/2}$ in $L^2(\mu,H)$.
By the weak lower semi-continuity of the norm, we get
		\begin{align*}
		\E\left[ \left| \nabla e^{-f/2} \right|^2  \right]&\leq \liminf \E\left[ \left| \nabla e^{-f_{nmkl}/2} \right|^2  \right]
		\end{align*}
		which implies that $\lim\limits_{n,l,k,m}\E\left[ \left| \nabla e^{-f_{nmkl}/2} \right|^2  \right]
		=\E\left[ \left| \nabla e^{-f/2} \right|^2  \right]$, and also $( \nabla e^{-f_{nmkl}/2} )$ converges  to $\nabla e^{-f/2}$ in $L^2(\mu,H)$ since $it$ converges weakly. Hence, 
		\begin{align*}
		\lim\limits_{n,l,k,m}\E\left[ |\nabla f_{nmkl}|^2 e^{-f_{nmkl}} \right]&=\lim\limits_{n,l,k,m}4\E\left[ \left| \nabla e^{-f_{nmkl}/2} \right|^2  \right]\\
		&=4\E\left[ \left| \nabla e^{-f/2} \right|^2  \right] \; =\; \E\left[ |\nabla f|^2 e^{-f} \right]
		\end{align*}
		and we get $\lim\limits_{n,l,k,m}E\left[ |\nabla f_{nmkl}|^2 e^{-f} \right]=\E\left[ |\nabla f|^2 e^{-f} \right].$ Moreover,  $(\nabla f_{nmkl})$ converges to $\nabla f$ in {$L^0(\rho)$}.  Hence, $(\nabla f_{nmkl})$ converges to  $\nabla f$ in $L^2(\rho).$
		Similar calculations show that $(\nabla f_{nmkl}\circ T_{nmkl})$ converges to $\nabla g\circ T$ in $L^2(\rho)$.
	\end{proof}
	
	We are ready to prove the regularity of backward Monge potential. In the next theorem, note that we do not assume Poincar\'{e} inequality since it is implied by $(1-c)$ convexity  \cite[Thm.6.2]{Convex_Feyel_Ustunel}.
	
	\begin{theorem}\label{4.1} 
		Let $(W,H,\mu)$ be an abstract Wiener space, $g\in\mathbb{D}_{2,1}$ and $f\in\mathbb{D}_{2,1}$ such that
		\begin{align*}
		\int_{W}|\nabla f|^2e^{-f}\;{d}\mu<\infty\;\text{and}\;\int_{W}|\nabla g|^2e^{-g}\;{d}\mu <\infty \nnumber{closable_condition2}
		\end{align*}
		and the function $f $ is $(1-c)$-convex function for some $c\in[0,1).$ Let $(\varphi, \psi)$ be  forward and backward  potentials to the Monge-Kantarovich problem
		with initial measure  and target measure 
		\begin{align*}
		\text{d}\rho=e^{-f}\text{d}\mu, \quad \quad \text{d}\nu=e^{-g}\text{d}\mu
		\end{align*}
		and  quadratic cost given by
		$$C(x,y)=\begin{cases} |x-y|^2_H &\mbox{if } x-y\in H \\
		\infty & \mbox{if }x-y\notin H .\end{cases}$$
		Then,   $\nabla^2\psi\in L^2(\nu,H\otimes H)$ and can be estimated by
		\begin{align*}
		\E_{\nu}\left[|\nabla^2\psi|_2^2 \right]
		&\leq \frac{3}{c}\left( \E_{\rho}\left[|\nabla\varphi|_H^2 \right] +\E_{\nu}\left[|\nabla g|_H^2\right]+\E_{\rho}\left[|\nabla f|_H^2\right]  \right)\, ,
		\end{align*}
		where $H\otimes H$ denotes the space of Hilbert-Schmidt operators on $H.$
	\end{theorem}
	\begin{proof}
		Thanks to condition \eqref{closable_condition2}, the Sobolev derivative is closable in $L^2(\rho)$ and $L^2(\nu)$. Define 
		\begin{align*}
		d\rho_{n}=e^{-f_{n}}\;d\mu, \quad \quad d\nu_{n}=e^{-g_{n} } d\mu
		\end{align*}
		with $f_n$ and $g_n$ as described before Lemma \ref{lemma5}.	The new probability measures $\rho_{n}$ and $\nu_{n}$ satisfy the sufficient conditions for   forward Monge potential to be smooth. Hence, if we apply Proposition \ref{regulary_smooth_bakcward}, we have		
		\begin{align*}
		c \E_{\nu_{n}}\left[|\nabla^2\psi_{n}|_2^2 \right]&\leq 3\left( \E_{\rho_{n}}\left[|\nabla\varphi_n|_H^2 \right] +\E_{\nu_{n}}\left[|\nabla g_{n} |_H^2\right]+\E_{\rho_{n}}\left[|\nabla f_{n}|_H^2\right]  \right)
		\end{align*}
		Lemmas \ref{lemma1}, \ref{lemma2}, \ref{lemma3} , \ref{lemma4} and \ref{lemma5} imply that the limits on right hand side exist and
		\begin{align*}
		\E_{\nu_{n}}\left[|\nabla^2\psi_{n}|_2^2 \right]&\leq \frac{3}{c}\left( \E_{\rho}\left[|\nabla\varphi|_H^2 \right] +\E_{\nu}\left[|\nabla g|_H^2\right]+\E_{\rho}\left[|\nabla f|_H^2\right]  \right).
		\end{align*}
		Those lemmas also imply that $(\sqrt{\theta_{{n}}} P_{1/n}\nabla\psi_n )$ converges weakly to $ \nabla^2 \psi .$  In view of the weak lower continuity of the norms, if we take a weak limit with respect to $n$, we obtain 
		\begin{align*}
		\E_{\nu}\left[|\nabla^2\psi|_2^2 \right]&\leq \liminf_{n} \E\left[|\sqrt{\theta_{{n}}} \nabla^2 P_{ {1}/{n} } \psi_{n} |_2^2 e^{-g} \right]\\
		&\leq \sup_{n}e^\frac{-2}{n} \E\left[\theta_{{n}} P_{ {1}/{n} }|\nabla^2\psi_{n}|_2^2 e^{-g_n}\right]\\
		&\leq \sup_{n}  \E\left[|\nabla^2\psi_{n}|_2^2 e^{-g_{n}}\right]\\
		&\leq \frac{3}{c}\left( \E_{\rho}\left[|\nabla\varphi|_H^2 \right] +\E_{\nu}\left[|\nabla g|_H^2\right]+\E_{\rho}\left[|\nabla f|_H^2\right]  \right)
		\end{align*}
		and the result follows. 
	\end{proof}		
	\begin{remark}
		The result of Lemma \ref{lemma5} is more than what is needed in the proof of Theorem \ref{4.1} Indeed, it holds that since $\E_{\nu_{n}}\left[|\nabla g_{n}|_H^2\right] \leq 8 \E_{\nu}\left[|\nabla g|_H^2\right]$ and  $\E_{\rho_{n}}\left[|\nabla f_{n}|_H^2\right] \leq 8\E_{\rho}\left[|\nabla f|_H^2\right]  $.
	\end{remark}
	
	\subsection{Monge-Amp\`ere Equation} \label{s4.3}
	We will prove that backward potential $\psi$ solves Monge-Amp\'ere equation. The idea of the proof is to start with finite dimension and take limits, as accomplished  through the following lemmas.
	
	\begin{lemma} Under assumptions of Theorem \ref{4.1},
		$(\mathcal{L}\psi_{n})$ converges to $\mathcal{L}\psi$ in the sense of distributions  and $\mathcal{L}\psi\in L^1(\nu)$.
	\end{lemma}
	\begin{proof}
		The proof of convergence in the sense of distributions follows from duality. Once we show that $(\mathcal{L}\psi_{n})$  is uniformly integrable, we are done. We will use an idea from  \cite{bogachev2013}. First we will show 
		\begin{align*}
		\sup_n \int \frac{(\mathcal{L}\psi_{n})^2}{1+|\nabla \psi_{n}|^2}\;d\rho< M<\infty
		\end{align*}
		then we will show the uniform integrability of $(\mathcal{L}\psi_{n})$ by using the convergence of $(\nabla\psi_{n})$ in $L^2(\nu).$ 
		Let $u$ be a decreasing function on $[0,+\infty].$ Observe that
		\begin{align*}
		\int \left( \mathcal{L}\psi_{n}\right) ^2 u(|\nabla \psi_{n}|^2) e^{-g}\;d\mu&=\int \left\langle \nabla\psi_{n} ,\nabla\left( \mathcal{L}\psi_{n} u(|\nabla \psi_{n}|^2) e^{-g}\right) \right\rangle _H\;d\mu
		\\
		&=\int \left\langle\nabla\psi_{n} ,\nabla\left( \mathcal{L}\psi_{n} u(|\nabla \psi_{n}|^2)  \right) \right\rangle_H e^{-g} \;d\mu
		\\
		&\qquad- \int \left\langle \nabla\psi_{n} ,\nabla g \right\rangle_H \mathcal{L}\psi_{n} u(|\nabla \psi_{n}|^2) e^{-g} \;d\mu
		\\
		&=\int \left\langle \nabla\psi_{n} ,\nabla\left( \mathcal{L}\psi_{n}  \right) \right\rangle _H u(|\nabla \psi_{n}|^2)  e^{-g} \;d\mu
		\\
		&\qquad+2\int \left\langle \nabla\psi_{n} ,\nabla^2 \psi(\nabla \psi) \right\rangle _H u'(|\nabla\psi_{n}|^2) \mathcal{L}\psi_{n} e^{-g} \;d\mu\\
		&\qquad- \int \left\langle \nabla\psi_{n} ,\nabla g \right\rangle_H \mathcal{L}\psi_{n} u(|\nabla \psi_{n}|^2) e^{-g} \;d\mu\\
		&=\int \left\langle \nabla\psi_{n} ,\mathcal{L}\nabla\psi_{n}  \right\rangle _H u(|\nabla \psi_{n}|^2)  e^{-g} \;d\mu
		\\
		&\qquad+\int \left\langle \nabla\psi_{n} ,\nabla\psi_{n}   \right\rangle_H u(|\nabla \psi_{n}|^2)  e^{-g} \;d\mu
		\\
		&\qquad+2\int \left\langle \nabla\psi_{n} ,\nabla^2 \psi(\nabla \psi) \right\rangle _H u'(|\nabla \psi_{n}|^2) \mathcal{L}\psi_{n} e^{-g} \;d\mu
		\\
		&\qquad- \int \left\langle \nabla\psi_{n} ,\nabla g \right\rangle_H \mathcal{L}\psi_{n} u(|\nabla \psi_{n}|^2) e^{-g} \;d\mu
		\end{align*}
		Finally, if we write the first line of the last term, we see that
		\begin{align*}
		\int \left( \mathcal{L}\psi_{n}\right) ^2 u(|\nabla \psi_{n}|^2) e^{-g}\;d\mu	&=\int |\nabla^2\psi_{n}  | _2^2  u(|\nabla \psi_{n}|^2) e^{-g} \;d\mu   \tag{{\Romannum{1}}} \\
		&\qquad-\int \left\langle \nabla^2\psi_{n}(\nabla\psi_{n}) ,\nabla g  \right\rangle_H u(|\nabla \psi_{n}|^2)  e^{-g} \;d\mu  \tag{\Romannum{2}}\\
		&\qquad+2\int | \nabla^2\psi_{n}(\nabla\psi_{n}) |_H^2 u'(|\nabla \psi_{n}|^2)  e^{-g} \;d\mu  \tag{\Romannum{3}}\\
		&\qquad+\int \left\langle \nabla\psi_{n} ,\nabla\psi_{n}   \right\rangle_H u(|\nabla \psi_{n}|^2)  e^{-g} \;d\mu \tag{\Romannum{4}}\\
		&\qquad+2\int \left\langle \nabla\psi_{n} ,\nabla^2 \psi(\nabla \psi) \right\rangle_H u'(|\nabla \psi_{n}|^2) \mathcal{L}\psi_{n} e^{-g} \;d\mu \tag{\Romannum{5}}\\
		&\qquad- \int \left\langle \nabla\psi_{n} ,\nabla g \right\rangle_H \mathcal{L}\psi_{n} u(|\nabla \psi_{n}|^2) e^{-g} \;d\mu \tag{\Romannum{6}}
		\end{align*}
		We know that $\Romannum{1},$ $\Romannum{4}$  are bounded and $\Romannum{3}$ is negative since $u$ is decreasing. We will show that the other integrals are also bounded when we choose $u$ properly. Observe that for all $\epsilon>0$ the Cauchy–Schwarz inequality implies 
		\begin{align*}
		|\Romannum{2}|&\leq \sqrt{\frac{1}{4\epsilon}\int |\nabla g|^2e^{-g}\;d\mu}\sqrt{4\epsilon \int |\nabla^2 \psi(\nabla \psi)|^2 u^2(|\nabla \psi_{n}|^2) e^{-g} \;d\mu} \\
		&\leq \frac{1}{4\epsilon}\int |\nabla g|^2e^{-g}\;d\mu +\epsilon \int |\nabla^2 \psi(\nabla \psi)|^2 u^2(|\nabla \psi_{n}|^2 e^{-g} \;d\mu
		\end{align*}
		\begin{align*}
		|\Romannum{6}|&\leq \sqrt{\frac{1}{4\epsilon}\int |\nabla g|^2e^{-g}\;d\mu}\sqrt{4\epsilon \int |\nabla \psi_{n} |^2(\mathcal{L}\psi_{n} )^2 u^2(|\nabla \psi_{n}|^2) e^{-g} \;d\mu} \\
		&\leq \frac{1}{4\epsilon}\int |\nabla g|^2e^{-g}\;d\mu +\epsilon \int |\nabla \psi_{n} |^2(\mathcal{L}\psi_{n} )^2 u^2(|\nabla \psi_{n}|^2) e^{-g} \;d\mu
		\end{align*}
		and
		\begin{align*}
		|\Romannum{5}|&\leq \sqrt{\frac{1}{\epsilon}\int (\mathcal{L}\psi_{n} )^2u(|\nabla \psi_{n}|^2 e^{-g}\;d\mu}\sqrt{4\epsilon \int \frac{(u')^2}{u}(|\nabla \psi_{n}|^2)|\nabla^2 \psi(\nabla \psi)|^2 |\nabla \psi|^2 e^{-g} \;d\mu} \\
		&\leq \epsilon \int (\mathcal{L}\psi_{n} )^2u(|\nabla \psi_{n}|^2) e^{-g}\;d\mu+ \frac{1}{\epsilon}\int \frac{(u')^2}{u}(|\nabla \psi_{n}|^2)|\nabla^2 \psi(\nabla \psi)|^2 |\nabla \psi|^2 e^{-g} \;d\mu.
		\end{align*}
		If we take $u(t)=\frac{1}{1+t}$ and $\epsilon=\frac{1}{4}$, we get
		\begin{align*}
		\frac{(u')^2}{u}(|\nabla \psi_{n}|^2)|)|\nabla \psi_{n}|^4\leq 1, \;\;\;\;\;   |\nabla \psi_{n}|^2 u(|\nabla \psi_{n}|^2)\leq 1
		\end{align*}
		and 
		\begin{align*}
		\epsilon u^2(|\nabla \psi_{n}|^2)+2u'(|\nabla \psi_{n}|^2)\leq 0
		\end{align*}
		Hence, we have \begin{align*}
		\sup_n \int \frac{(\mathcal{L}\psi_{n})^2}{1+|\nabla \psi_{n}|^2}\;d\rho< M<\infty
		\end{align*}
		where 
		\begin{align*}
		M= \int |\nabla g|^2e^{-g}\;d\mu + 2\sup_n \int |\nabla \psi_{n}|^2 e^{-g}\;d\mu +10\sup_n \int |\nabla^2 \psi_{n}|^2_2 e^{-g}\;d\mu<\infty
		\end{align*}
		Note that $(\nabla \psi_{n})$  is uniformly integrable with respect to $\nu.$ So for every $\epsilon>0$ there exists $\delta >0$ such that 
		$\rho(E)<\delta$ implies  $\E_\nu [|\nabla \psi|^21 _E]<\frac{\epsilon^2}{1+4M}.$
		Define
		\begin{align*}
		1_{E_1}=\left\lbrace |\mathcal{L}\psi_{n}|\leq \frac{\epsilon}{2M}\frac{|\mathcal{L}|\psi_{n}|^2}{1+|\nabla \psi_n|^2}\right\rbrace \\
		1_{E_2}=\left\lbrace |\mathcal{L}\psi_{n}|> \frac{\epsilon}{2M}\frac{|\mathcal{L}|\psi_{n}|^2}{1+|\nabla \psi_n|^2}\right\rbrace 
		\end{align*}
		We have $\E_\nu[|\mathcal{L}\psi_{n}|1_E]=\E_\nu[|\mathcal{L}\psi_{n}|1_{E_1}]+\E_\nu[|\mathcal{L}\psi_{n}|1_{E_2}]<\frac{\epsilon}{2}+\frac{\epsilon}{2}.$ Hence  $(\mathcal{L}\psi_{n})$ is uniformly integrable with respect to $\nu$. Therefore, $(\mathcal{L}\psi_{n})$ converges weakly to $\mathcal{L}\psi$ in $L^1(\nu).$
	\end{proof}

	\begin{lemma}\label{lemma9} Under assumptions of Theorem \ref{4.1},
		$(f_n\circ S_n)$ converges to $f\circ S$ in $L^1(\nu)$, where $S=T^{-1}=I_W+\nabla \psi$.
	\end{lemma}
	\begin{proof}We have
		\begin{align*}
		\int |f_n\circ S_n-f\circ S|\;d\nu\leq \int |f_n\circ S_n-f\circ S_n|_H\;d\nu + \int |f\circ S_n-f\circ S|_H\;d\nu 
		\end{align*}
		By Fatou's lemma 
		\begin{align*}
		\lim	\int |f_n\circ S_n-f\circ S_n|\;d\nu \leq  \lim \int |f_n-f|e^{-f_n}\;d\mu =0
		\end{align*}
		Let $\hat{f}\in C_b(W)$ such that $|f-\hat{f}|_{L^1(\rho)}<\epsilon.$ Then, again using Fatou's lemma
		\begin{align*}
		\lim \int |f\circ S_n-f\circ S|\;d\nu &\leq \lim \int |f\circ S_n-\hat{f}\circ S_n|e^{-g_n}\;d\mu+\lim \int |\hat{f}\circ S_n-\hat{f}\circ S|e^{-g}\;d\mu \\ &\qquad+   \int |\hat{f}\circ S-f\circ S|e^{-g}\;d\mu\\
		&\leq\lim \int |f-\hat{f}|e^{-f_n}\;d\mu+\lim\int |\hat{f}\circ S_n-\hat{f}\circ S| e^{-g}\;d\mu \\ &\qquad+   \int |\hat{f}-f|e^{-f}\;d\mu\\
		&\leq2|f-\hat{f}|_{L^1(\rho)}< 2\epsilon.
		\end{align*}
		Therefore, $\lim \int |f_n\circ S_n-f\circ S|\;d\rho =0$.
	\end{proof}
	\begin{lemma}\label{lemma8} 
		Under assumptions of Theorem \ref{4.1},
		\begin{align*}
		\lim\limits_{n} E_\nu \left[ |\nabla^2\psi-\nabla^2\psi_{n}|_2 \right]=0 
		\end{align*}
		where $|\cdot|_2$ denotes the  Hilbert-Schmidt operator norm.
	\end{lemma}
	\begin{proof}
		Let $m,$ $n$ be two integers such that $m>n$ and $\pi^n_m$ be the orthogonal projection from $H_m$ onto $H_n.$  Then  \cite[Thm.2.5]{fang} implies 
		\begin{align*}
		|\nabla^2 (\psi_{n}\circ\pi^m_n )-\nabla^2 \psi_{m}|_{L^1(\rho_{m},H_m\otimes H_m)}^2 \leq & \; C_1 \int_{H^m }(f_m-f_n\circ \pi^n_m)\; d \rho_m \\
		&+\frac{C_2}{\epsilon}\int_{H_m} |\nabla g_m- \nabla(g_n\circ \pi^n_m)|^2\;d\nu_m  
		\end{align*}
		for some constants $C_1, C_2 > 0$. 
		If we take the limit with respect to $m$ first, and then with respect to $n$, we get the result.
	\end{proof}
	
	\begin{theorem} \label{4.2}
		Let $(W,H,\mu)$ be an abstract Wiener space, $f\in\mathbb{D}_{2,1}$ and $g\in\mathbb{D}_{2,1}$ such that
		\begin{align*}
		\int_{\mathbb{R}^d}|\nabla f|^2e^{-f}\;{d}\mu<\infty\;\text{and}\;\int_{\mathbb{R}^d}|\nabla g|^2e^{-g} \;d\mu <\infty
		\end{align*}
		the function $f $ is $(1-c)$-convex function for some $c\in[0,1).$  
		Let $\psi$ be  backward  potential of Monge-Kantarovich problem
		with initial   and target measures
		\begin{align*}
		\text{d}\rho =e^{-f}\text{d}\mu , \quad \quad
		\text{d}\nu  =e^{-g}\text{d}\mu
		\end{align*}
		and  quadratic cost. Then, the backward  potential $\psi$ solves Monge-Amp\`ere equation
		\begin{align*}
		e^{-g}=e^{-f\circ S}{\det}_2\left( I_{H} +\nabla^2\psi\right) \exp \left[ -\mathcal{L}\psi -\frac{1}{2}|\nabla\psi|^2_H \right] 
		\end{align*}
		$\nu$-a.s., where  $S=T^{-1}=I_W+\nabla \psi$.
	\end{theorem}
	\begin{proof}
		By the finite dimensional result in \cite{villani2003topics}, we have 
		\begin{align*}
		e^{-g_n}=e^{-f_n\circ S_n}{\det}_2\left( I_{H} +\nabla^2\psi_n\right) \exp \left[ -\mathcal{L}\psi_n -\frac{1}{2}|\nabla\psi_n|^2_H \right] .
		\end{align*}
		 We already know that $e^{-g_n}\to e^{-g}$. From approximation lemmas in Subsection \ref{s4.1}, as $n\to \infty$, $   \sqrt{\theta_{n}} P_{1/n}\nabla\psi_n\to \nabla \psi  $ $\nu$-a.s.  Due to the continuity of $P_t$ in $t$ \cite[Prop.2.1]{SPS_1995} and since $\theta_n \rightarrow 1$, it follows that $ \nabla\psi_n\to \nabla \psi  $ $\nu$-a.s. By lemmas   \ref{lemma9} and \ref{lemma8}, we also know that $ e^{-f_n\circ T_n}\to e^{-f\circ T}, $   and  $\nabla^2\psi_{n}\to\nabla^2 \psi$ and the sequence  $(-\mathcal{L}\psi_n)$ has a limit, say $A$, $\nu$-a.s. Therefore, we have  
		\begin{align*}
		e^{-g}=e^{-f\circ S}{\det}_2\left( I_{H} +\nabla^2\psi\right) \exp \left[ A-\frac{1}{2}|\nabla\psi|^2_H \right] .
		\end{align*}
		Once we show that $-\mathcal{L}\psi=A$ holds $\nu$-a.s., we are done. Indeed, since $(-\mathcal{L}\psi_n)$ is uniformly integrable, for every cylindrical function $\xi$, we have
		\begin{align*}
		\int A\xi e^{-g}\;d\mu&=\lim_n \int -\mathcal{L}\psi_n\xi e^{-g}\;d\mu\\
		&=\lim_n \int-\langle\nabla \psi_{n}, \nabla\xi-\xi\nabla g\rangle_H\psi_n e^{-g}\;d\mu\\
		&=\lim_n \int-\langle\nabla \psi, \nabla\xi-\xi\nabla g\rangle_H\psi_n e^{-g}\;d\mu\\
		&= \int -\mathcal{L}\psi\xi e^{-g}\;d\mu
		\end{align*}
		so $-\mathcal{L}\psi=A$ holds $\nu$-a.s. 
	\end{proof}

\noindent{{\bf Acknowledgements.} This work is supported by TUBITAK Project 118F403. The authors are grateful to A. Suleyman Ustunel for his helpful guidance and instructive comments on the manuscript. 
	
	\bibliographystyle{plain}
	\bibliography{monge.bib}
\end{document}